\documentclass[12pt]{article}

\usepackage{amsfonts}
\usepackage{graphics}
\usepackage{amsmath}
\usepackage{amscd}
\usepackage{amssymb}
\usepackage{rlepsf}
\usepackage{color}

\usepackage{euscript}

\newtheorem{Theorem}{Theorem}[section]

\newtheorem{Corollary}[Theorem]{Corollary}
\newtheorem{Example}[Theorem]{Example}
\newtheorem{Definition}[Theorem]{Definition}

\newtheorem{Remark}[Theorem]{Remark}

\newtheorem{Lemma}[Theorem]{Lemma}
\newtheorem{Proposition}[Theorem]{Proposition}

\newtheorem{Fundamental Theorem}{Fundamental Theorem}

\newenvironment{Proof}[1][Proof]{\textbf{#1.} }{\ \rule{0.5em}{0.5em}}

\def \A {\mathcal{A}}
\def \B {\mathcal{B}}

\def \b {\beta}

\def\leq{\leqslant} 
 
\def \ge {\geqslant}
\def \bc {{\bf c}}
\def \C {\mathcal{C}}
\def \CRS {\mathrm{CRS}}

\def \D {\Delta}
\def \d {\partial}
\def \f {\phi}
\def \G {\mathcal{G}}

\def \id {\mathrm{id}}
\def \l {\lambda}
\def \L {\Lambda}
\def \N {\mathbb{N}}
\def \n {\eta}
\def \p {\psi}
\def \ra {\xrightarrow}

\def \s {\scriptstyle}
\def \S {\mathcal {S}}
\def \SIMP {\mathrm{SIMP}}
\def \t {\triangleright}
\def \tn {\otimes}
\def \TOP {\mathrm{TOP}}
\def \w {\omega}
\def \Z {\mathbb{Z}}
\begin{document}

\title{On Yetter's Invariant and an Extension of the Dijkgraaf-Witten Invariant to Categorical Groups}

\author{ Jo\~{a}o  Faria Martins\footnote{Financed by  Funda\c{c}\~{a}o para a Ci\^{e}ncia e Tecnologia (Portugal),
post-doctoral grant number SFRH/BPD/17552/2004, part of the research project
POCTI/MAT/60352/2004 (``Quantum Topology''), also financed by FCT. } \footnote{Also: Universidade Lus\'{o}fona de Humanidades e Tecnologia, Av do Campo Grande, 376, 1749-024, Lisboa, Portugal. }\\ \footnotesize\it  {Departamento de Matem\'{a}tica, Instituto Superior T\'{e}cnico,}\\ {\footnotesize\it Av. Rovisco Pais, 1049-001 Lisboa, Portugal.}\\  {\footnotesize\it jmartins@math.ist.utl.pt}\\  \\Timothy Porter 
\\  \footnotesize\it {Department of Mathematics,} \footnotesize\it { 
                               University of Wales, Bangor,} 
\\ \footnotesize\it {
                               Dean St.,
                               Bangor,
                               Gwynedd LL57 1UT,
                               UK.}\\  
      \footnotesize\it {   t.porter@bangor.ac.uk}}

\date{\today}

\maketitle

\begin{abstract}
We give an interpretation of Yetter's Invariant of manifolds $M$ in terms of
the homotopy type of the function space $\TOP(M,B(\G))$, where $\G$ is a
crossed module and $B(\G)$ is  its classifying space.  From this formulation,  there follows  that Yetter's invariant depends only on the homotopy type of $M$, {and the weak  homotopy type of the crossed module $\G$.} We  use this
interpretation  to define  a twisting of Yetter's Invariant by cohomology
classes of crossed modules,  defined as cohomology classes of their
classifying spaces,  in the form of a state sum invariant. In particular, we
obtain an extension of the  Dijkgraaf-Witten Invariant of manifolds to categorical groups. The straightforward extension to crossed complexes is also considered. 
\end{abstract}
\begin{center}
\end{center}

\section*{Introduction}

Let $G$ be a finite group. In the context of Topological Gauge Field Theory, R. Dijkgraaf
and E. Witten defined in \cite{DW} a 3-dimensional manifold invariant for each
3-dimensional cohomology class of $G$. When $M$ is a  triangulated manifold, the Dijkgraaf-Witten Invariant can be expressed in terms of a state sum model, reminiscent of Lattice Gauge Field Theory. State sum invariants of manifolds were popularised  by V.G. Turaev and
O.Ya. Viro, because of  their construction of a non-trivial closed 3-manifold state sum invariant
from quantum $6j$-symbols (see \cite{TV}),  regularising the
divergences of the celebrated Ponzano-Regge Model (see \cite{PR}), which is  powerful
enough to distinguish between manifolds which are homotopic but not {homeomorphic}. For categorically inclined introductions to state sum invariants of manifolds we refer the reader to \cite{T,BW,Mk1}.

 State sum models, also known as spin foam models, appear frequently in the
 context of 3- and 4-dimensional Quantum Gravity, as well as Chern-Simons
Theory, making them particularly interesting. See for example \cite{B,Bae} for reviews.

Let us recall the construction of the Dijkgraaf-Witten Invariant. 
Let $G$ be a finite group.  A 3-dimensional group cocycle is a map $\w\colon    G^3 \to
U(1)$ verifying the cocycle condition: \[\w(X,Y,Z)\w(XY,Z,W)^{-1} \w(X,YZ,W) \w(X,Y,ZW)^{-1}
\w(Y,Z,W)=1,\] for any $ X,Y,Z,W \in G$. These group 3-cocycles can  be
seen as 3-cocycles of the classifying space $B(G)$ of $G$, by considering the
usual simplicial structure on $B(G)$; see for example \cite{We}, Example $8.1.7$.

\begin{figure}
\centerline{\relabelbox 
\epsfysize 2cm
\epsfbox{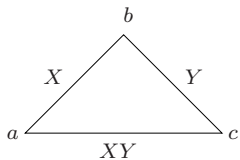}
\relabel {X}{$\s{X}$}
\relabel {Y}{$\s{Y}$}
\relabel {Z}{$\s{XY}$}
\relabel {1}{$\s{a}$} 
\relabel {2}{$\s{b}$}
\relabel {3}{$\s{c}$}
\endrelabelbox }
\caption{\label{simplex1} Flatness condition imposed on $G$-colourings of a 
  triangulated manifold $M$. Here $X,Y \in G$ and $a<b<c$ are vertices of $M$. }
\end{figure}

\begin{figure}
\centerline{\relabelbox 
\epsfysize 2cm
\epsfbox{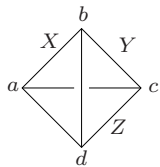}
\relabel {X}{$\s{X}$}
\relabel {Y}{$\s{Y}$}
\relabel {Z}{$\s{Z}$}
\relabel {A}{$\s{a}$} 
\relabel {B}{$\s{b}$}
\relabel {C}{$\s{c}$}
\relabel {D}{$\s{d}$}
\endrelabelbox }
\caption{\label{action} Edges of a tetrahedron of $M$ whose colours determine
  the colouring of all other edges. Here $a<b<c<d$ are  vertices of $M$.}
\end{figure}

 Let $M$ be an  oriented closed connected piecewise linear manifold. Choose a triangulation of
 $M$, and consider a base point $*$ for $M$ so that $*$ is a vertex of $M$.  Suppose that we are provided with a total order on the set of vertices
 of $M$. A $G$-colouring of $M$ is given by an assignment of an element of $G$ to each
 edge of $M$, satisfying the flatness condition shown in figure
 \ref{simplex1}. A $G$-colouring of $M$ therefore defines a morphism $\pi_1(M,*)
 \to G$. Let $T=(abcd)$, where $a<b<c<d$, be a tetrahedron of $M$. The
 colouring of each edge  of $T$ is determined by the colouring of the edges
 $(ab),(bc)$ and $(dc)$, of $T$; see figure \ref{action}.   Let $\w{\colon}G^3 \to U(1)$ be a
 3-dimensional group cocycle, representing some 3-dimensional cohomology class  of $G$.
The ``topological action'' $S(T,\w)$ on a tetrahedron $T$ is $\w(X,Y,Z)^{\pm 1}$, where
 $X,Y,Z\in G$ are as in   figure \ref{action}. Here the exponent on $\w(X,Y,Z)$ is $1$ or
 $-1$ according to whether the orientation on $T$ induced by the total order
 on the  set of vertices of $M$ coincides or not with the orientation of $M$.
 The expression of  the Dijkgraaf-Witten Invariant is the state sum:
\[DW(M,\w)=\frac{1}{\#G^{n_0}}\sum_{G\textrm{-colourings}} \quad 
\prod_{\textrm{tetrahedra } T  } S(T,\w),\]
where $n_0$ is the number of vertices of the chosen triangulation of $M$. This state sum is triangulation independent. 

A  homotopy theoretic expression for the Dijkgraaf-Witten Invariant is:
\[DW(M,\w)={\frac{1}{\# G}}\sum_{f \in [(M,*),(B(G),*)]} \left <o_M,f^*(\w) \right>,\]
 where $o_M \in H_3(M)$ denotes the orientation class of $M$. The notation
 $[(M,*),(B(G),*)]=\pi_0\big (\TOP((M,*),(B(G),*)) \big)$ stands for the set of based homotopy classes of maps $(M,*)  \to (B(G),*)$, where, as before, $B(G)$ denotes the classifying space of $G$. The equivalence between the two formulations is a consequence of
 the fact that there exists a one-to-one correspondence between group morphisms
 $\pi_1(M,*) \to G$ and pointed  homotopy classes of maps $(M,*) \to (B(G),*)$.

In \cite{Y1}, D.N. Yetter defined a 3-dimensional TQFT, which included  the
$\w=1$ case of the Dijkgraaf-Witten Invariant, though he did not incorporate group cohomology classes into his TQFT. This construction
appeared upgraded  in \cite{Y2}, where D.N. Yetter extended his framework to handle  categorical
groups. Recall that categorical groups are equivalent to crossed modules. 
 They are often easier to consider when handling the objects at a theoretical
 level, although for calculation the crossed modules have advantages.  We will
 often pass from one description to the other, usually without comment.  The
 important fact is that they both are algebraic models of homotopy 2-types and
 so generalise groups, which model homotopy 1-types.

Let $M$ be a piecewise linear compact manifold equipped with a triangulation. The definition  of Yetter's Invariant (which we will formulate carefully in this article)
is analogous to the definition of the Dijkgraaf-Witten Invariant for the case  $\w=1$. However, in Yetter's construction,  in addition to colourings of the edges of  $M$, now by
objects of a categorical group,  we also have colourings on
the faces the triangulation of $M$ by morphisms of the chosen categorical group, and the flatness condition is transported to the
tetrahedrons of $M$.

The meaning of Yetter's construction was elucidated by the second author in \cite{P1,P2},  where an extension  of Yetter's Invariant to handle general
models of $n$-types appears. The construction of M. Mackaay in \cite{Mk2} is related, conjecturally,
with the $n=3$ case of the  construction in \cite{P2}, with a further twisting by
cohomology classes of $3$-types, in the Dijkgraaf-Witten way.

In \cite{FM1,FM2}, D.N. Yetter's ideas were applied to 2-dimensional knot theory,
and yielded interesting invariants of knotted embedded  surfaces in $S^4$, with nice calculational properties. In \cite{FM2,FM3} the
general formulation  of Yetter's Invariant for CW-complexes and crossed complexes appeared.

One of the open problems posed by D.N. Yetter at the end of \cite{Y2} was whether his
construction could, similarly to the  Dijkgraaf-Witten Invariant, be twisted by
cohomology classes, and he asked for the right setting for doing this. In this article we formulate such  a twisting of
Yetter's Invariant, thereby extending the Dijkgraaf-Witten Invariant of 3-manifolds to categorical groups.
Note that the (co)homology of a crossed module is defined as being the (co)homology of its classifying space; see for example \cite{E,Pa}.

In the course of this work, we  prove a new homotopy theoretical interpretation of
Yetter's Invariant (Theorem \ref{main1}), therefore giving a new proof of its existence; see \cite{P1,P2} for alternative
interpretations. Another important consequence of this formulation is the homotopy invariance of Yetter's Invariant, as well as {the fact that Yetter's Invariant depends only on the weak homotopy type of the chosen crossed module.}

Our interpretation of Yetter's Invariant
extends in the obvious way to the homologically twisted case (Theorem \ref{main2}), and has
close similarities with the interpretation of the  Dijkgraaf-Witten Invariant shown above. The main difference is that it is set in the unbased case, and depends on the homotopy groups of each connected component of the function space $\TOP(M,B(\G))$. Here $B(\G)$ is the classifying space of the crossed module $\G$. 
  Our construction does not depend on  the dimension of the  manifolds, and can also be adapted  to the general crossed complex context, which we  expect to be closely related with the invariants appearing in \cite{P2,Mk2}.   This extension of Yetter's  invariant to crossed complexes, as well as its homotopy invariance and geometric interpretation appears in this article; {see \ref{ABC}.}

As we mentioned before,  the category of categorical groups is equivalent to the
category of crossed modules, which are particular cases of crossed
complexes. Crossed complexes $\A$ admit classifying spaces $B(\A)$, which
generalise  Eilenberg-Mac Lane spaces. 
  The classifying spaces of crossed complexes $\A$ were  extensively studied by R. Brown and P.J. Higgins. One of their main results is the total description in algebraic terms of the weak homotopy type of the function space
$\TOP(M,B(\A))$, where $M$ is a CW-complex.  This result together with the General van  Kampen Theorem (also due to R. Brown and P.J. Higgins) are the main tools that we will use in this article. 

We will make our exposition as self contained as possible, giving brief
description of the concepts of the theory of  crossed complexes that we will
need, constructions which are due mainly to R. Brown and P.J. Higgins.

This article should be compared to \cite{P3,PT}, where similar ideas are used to define Formal Homotopy Quantum Field Theories with background a 2-type.

\tableofcontents

\section{Crossed Modules and Yetter's Invariant}

\subsection{Crossed Modules and Crossed Complexes}
\subsubsection{Definition of Crossed Modules}
Let $G$ be a groupoid with object set $P$. We denote the source and target
maps of $G$ by $s,t\colon   G \to P$, respectively. If a groupoid $E$  is  totally disconnected, we denote both maps (which are by definition equal) by $\b\colon   E \to P$.
\begin{Definition}
Let $G$ and $E$ be groupoids, over the same set $P$, with $E$ totally
disconnected.
A  crossed module $\G=(G,E,\d,\t)$ is given by an object
preserving  groupoid morphism $\d\colon    E\to G$ and a left  groupoid action $\t$ of $G$ on $E$  by automorphisms. The conditions on $\t$ and $\d$ are:
\begin{enumerate}
\item $\d(X \t e)=X\d(e)X^{-1}$, if $X \in G$ and $e \in E$ are such that $t(X)=\b(e)$,
\item $\d(e) \t f=e f e^{-1}$, if $ e, f \in E$ verify $\b(e)=\b(f)$.
\end{enumerate}
\end{Definition}
See \cite{BI}, section 1.
Notice that for any $e \in E$ we must have $s(\d(e))=t(\d(e))=\b(e)$.

A crossed module is called {\it reduced} if $P$ is a singleton. This implies that both $G$ and $E$ are groups.  A morphism $F=(\f,\psi)$ between the crossed modules $\G=(G,E,\d,\t)$ and
$\G'=(G',E',\d',\t')$ is given by a pair of groupoid morphisms $\f\colon   G\to G'$ and
$\psi\colon    E \to E'$, making the diagram:
\begin{equation*}
\begin{CD}
E @>\p>> E' \\
@V\d VV  @VV\d' V\\
 G @>>\f > G'
\end{CD}
\end{equation*}
commutative. 
In addition we must have:   \[
t(X)=\b(e) \implies \f(X)\t' \psi (e)=\psi(X \t e), \forall X \in G, \forall e \in E.\]

There exists an extensive literature on (reduced) crossed modules. We refer for example to 
\cite{BA1,BA2,BHu,FM2}. The non-reduced case, important for this article, is considered in \cite{B1,B4,BH1,BH2,BSi,BI}. 

\subsubsection{Reduced Crossed Modules and Categorical Groups}

Let  $\G=(G,E,\d,\t)$ be a reduced crossed module. Here $\d \colon    E \to G$ is a
group morphism,  and
$\t$ is a left action of  $G$ on $E$  by automorphisms.  We can define a strict monoidal category $\C(\G)$ from $\G$. The set of objects of $\C(\G)$ is given by all elements of $G$. Given an $X \in G$, the set of all morphisms with source $X$ is in bijective correspondence with  $E$, and the target of $e\in E$ is $\d(e^{-1}) X$. In other words a morphism in $C(\G)$   ``looks like'' 
 $X \ra{e} \d(e^{-1}) X$. Given $X \in G$ and $e,f \in E$ the composition 
 
\[X \ra{e} \d(e^{-1})X \ra{f} \d(f^{-1}) \d(e^{-1}) X\] is:
 \[X \ra{ef} \d(ef)^{-1}X.\]   

The tensor product has the form:
\begin{equation*}
\begin{CD}
\begin{CD} X \\ @VVeV \\\d(e^{-1})X \end{CD} \tn \begin{CD} Y \\ @VVfV
  \\\d(f^{-1})Y \end{CD} =\begin{CD} XY \\ @VV (X\t f)e V \\\d(e^{-1})X\d(f^{-1})Y \end{CD}
\end{CD}.
\end{equation*}
From the definition of a crossed module, it is easy to see that we have indeed
defined a strict  monoidal category. The
tensor category $\C(\G)$ is a categorical group (see \cite{BM} and
\cite{BSi}). It is well known that the categories of crossed modules and of
categorical  groups are equivalent (see \cite{BSi,BSp,BM}).  This construction is an old one. We skip further details on
this connection since the  point we want to emphasise is that  we can
construct a strict monoidal category $\C(\G)$ which is moreover a categorical group  from any crossed module $\G$.  This also justifies
the title of this article. Although we will define Yetter's  invariant using the crossed module
form of the data, it is, we feel, more natural to think of the categorical
group as being the coefficient object.

\subsubsection{Crossed Complexes}

A natural generalisation  of the concept of  a crossed module is a crossed complex.

\begin{Definition}
A crossed complex $\A$ is given by a complex of groupoids over the same set
$A_0$ (in other words $A_0$ is the object set of all groupoids), say:
\[...\ra{\d_{n+1}} A_n \ra{\d_n}  A_{n-1} \ra{\d_{n-1}}A_{n-2}
\ra{\d_{n-2}}...\ra{\d_3} A_2 \ra{\d_2} A_1 ,\]
where all boundary maps $\d_n, n \in \N$ are object preserving, and the  groupoids
$A_n, n\ge 2$ are   totally disconnected. The remaining conditions are:
\begin{enumerate}

\item  For each $ n >1$, there exists a left groupoid action $\t=\t_n$ of the groupoid $A_1$ on $A_n$, by automorphisms, and all the boundary maps $\d_n, n=2,3,...$ are $A_1$-module morphisms.

\item The groupoid map  $A_2 \ra{\d_2} A_1$ together with the action $\t$ of
  $A_1$ in $A_2$ defines a crossed module. 
\item The groupoid $A_n$ is abelian if $n>2$.

\item $\d_2(A_2)$ acts trivially on $A_n$ for $n>2$.
\end{enumerate}
For a positive integer $L$, an $L$-truncated crossed complex is a crossed complex such that $A_n=\{0\}$ if $n>L$. In particular crossed modules are equivalent to $2$-truncated crossed complexes.  A crossed complex is called reduced if $A_0$ is a singleton. 
\end{Definition}
 Morphism of crossed complexes are defined in a a similar way to morphisms of crossed modules.

A natural example of a crossed complex is the following  one, introduced by
A.L. Blakers in \cite{Bl}, in the reduced case, and extensively studied by J.H.C.
Whitehead in \cite{W4,W5}, and H.J. Baues, see \cite{BA1,BA2}. The non-reduced case, important for this work, is treated in \cite{B1,B4,BH2,BH3,BH4,BH5}, {which we refer to for further details.} 
\begin{Example}
Let $M$ be a locally path connected space, so that each connected component of $M$ is path connected. A filtration of $M$ is an
increasing sequence $\{M_k,k=0,1,2,...\}$ of locally path-connected subspaces
of $M$, whose union is $M$. In addition, we suppose that 
$M_0$ has a non-empty intersection with each connected component of $M_k$, for each $k\in\N$.  Then the sequence of groupoids $\pi_n(M_n,M_{n-1},M_0)$, with the
obvious boundary maps, and  natural left actions of the fundamental groupoid{\footnote{ {This should not  be confused with the relative homotopy set $\pi_1(M_1,M_0,*)$.}}} $\pi_1(M_1,M_0)$, {with set of base points $M_0$},  on them
is a crossed complex, which we denote by $\Pi(M)$, and call the ``Fundamental
Crossed Complex of the Filtered Space $M$''. Note that the assignment $M \mapsto
\Pi(M)$, is functorial.
\end{Example}

\begin{Remark}
Let $M$ be a CW-complex. The notation $\Pi(M)$ will always mean the fundamental crossed complex
of the skeletal filtration $\{M^k, k\in\N\}$ of $M$, for which $M^k$ is the
$k$-skeleton of $M$. Note that it makes sense to consider the fundamental
crossed complex of this type of  filtrations.  
\end{Remark}

It is easy to show  that crossed complexes and their morphisms
form a category. 
We will usually denote a crossed complex $\A$ by $\A=(A_n,\d_n,\t_n)$, or more
simply by $(A_n,\d_n)$, or  even $(A_n)$. A morphism $f\colon    \A \to \B$ of
crossed complexes will normally be denoted by $f=(f_n)$. 

Let $\A=(A_n,\d_n,\t_n)$ be a crossed complex. For any $a \in A_0$ and any $n \ge 1$, denote
$A_n(a)$ as being  $\{a_n \in A_n\colon    s(a_n) =a\}$, which is therefore a group if
$n \ge 2$. Denote also $A_1(a,b)=\{a_1 \in A_1\colon   s(a_1)=a, t(a_1)=b\}$, where
$a,b \in A_0$.

\subsubsection{The General van Kampen Theorem}\label{general}
 The category of crossed complexes  is a category with colimits; see
 \cite{BH2}, section $6$.  
Under mild conditions, the functor $\Pi$ from the category of filtered
 spaces to the category of crossed complexes preserves colimits. This fact is known as the ``General van Kampen Theorem'', and is due to
 R. Brown and P.J. Higgins; see \cite{B1,B3,B4,BH1,BH2,BH3,BH4}, and also \cite{BA1,BA2}. A useful form of the
 General van Kampen Theorem is the following:
\begin{Theorem}  {\bf (R. Brown and P.J. Higgins)}
Let $M$ be a CW-complex. Let also $\{M_\l,\l \in \L\}$ be a covering of $M$ by
subcomplexes of $M$. Then the following diagram:
\[\bigsqcup_{(\l_1,\l_2) \in \L^2} \Pi(M_{\l_1} \cap M_{\l_2}) {\substack{ \longrightarrow\\ \longrightarrow}} \bigsqcup_{\l
  \in \L}\Pi( M_\l) \longrightarrow \Pi(M), 
\]
 with arrows induced by inclusions, is a coequaliser diagram in the category of crossed complexes.
\end{Theorem} 
This is corollary $5.2$ of \cite{BH2}.

\subsubsection{Homotopy of Maps from Crossed Complexes into Crossed Modules}\label{HomMaps}
The main references here are \cite{BI,BH4}; see also  \cite{B1,B4,W4}. Note that our conventions are slightly different. 
Let $G$ and $G'$ be groupoids. Let also $E'$ be a  totally disconnected groupoid
on which $G'$ has a left action. Let $\f\colon   G \to G'$ be a groupoid morphism. A
map $d\colon   G \to E'$ is said to be a $\f$-derivation if: 
\begin{enumerate}
\item $\b(d(X))=t(\f(X)),\forall X \in G$.
\item $d(XY)=\big (\f(Y)^{-1} \t d(X)\big) d(Y)$ if $X,Y \in G$ are such that $t(X)=s(Y)$. 
\end{enumerate}

\begin{Definition}
Let $\A=(A_n)$ be a crossed complex and $\G=(G,E,\d,\t)$ be a crossed module. Let also $f\colon   \A \to \G$ be a morphism of crossed complexes. An $f$-homotopy $H$ is given
by two maps: 
\[H_0\colon   A_0 \to G, \quad H_1\colon   A_1 \to E,\]
such that:
\begin{enumerate}
\item $s(H_0(a))=f_0(a), \forall a \in A_0$,
\item $H_1\colon   A_1 \to E$ is an $f_1$-derivation, here $f=(f_n)$,
\end{enumerate}
in which case we put $f=s(H)$.
\end{Definition}
\begin{Proposition}
Under  the conditions of the previous definition,  suppose that  $H$ is an $f$-homotopy. The map $g=t(H)\colon   \A \to \G$ such that:
\begin{enumerate}
\item $g_0(a)=t(H_0(a)), \forall a \in A_0,$
\item $g_1(X)=H_0(s(X))^{-1}f_1(X)(\d \circ H_1)(X)H_0(t(X)), \forall X \in A_1, $
\item $g_2(e)=H_0(\b(e))^{-1} \t \left (f_2(e)(H_1\circ \d)(e)\right), \forall e \in A_2$,
\end{enumerate} 
is a morphism of crossed complexes. Moreover, this construction defines a
groupoid $\CRS_1(\A,\G)$, whose  objects  are given by all morphisms $\A \to \G$, and given two morphisms $f,g\colon   \A \to \G$, a morphism $H\colon   f \to g$ is an $f$-homotopy $H$ with $t(H)=g$.
\end{Proposition}
For each $f \in \CRS_0(\A,\G)$ (the set of crossed complex  morphisms $\A \to \G$), denote the
set of all $f$-homotopies by $\CRS_1(\A,\G)(f)$. Analogously, if $f,g\in \CRS_0(\A,\G)$, we denote $\CRS_1(\A,\G)(f,g)$ as being the set of elements of the groupoid $\CRS_1(\A,\G)$ with source $f$ and target $g$.
\begin{Definition}
Let $\A=(A_n)$ be  a crossed complex and $\G=(G,E,\d,\t)$ be a crossed module. Let
also $f\colon   \A \to \G$ be a morphism. A 2-fold $f$-homotopy is a map $I\colon   A_0 \to E$
such that $\b(I(a))=f_0(a), \forall a \in A_0$. As usual $f=(f_n)$. 
\end{Definition}
Denote  $\CRS_2(\A,\G)(f)$ as being the set of all 2-fold $f$-homotopies. It is a group with pointwise multiplication as product. Therefore the disjoint union:
\[\CRS_2(\A,\G)\doteq \bigsqcup_{f \in \CRS_0(\A,\G)} \CRS_2(\A,\G)(f)\]
is a totally disconnected groupoid. 
\begin{Proposition}
Under  the conditions of the previous definition,  let $I$ be a 2-fold $f$-homotopy. Then  $\d I=(H_1,H_2)$, where:
\begin{enumerate}
\item $H_0(c)=\d(I(c)), \forall c \in A_0,$
\item $H_1(X)=f_1(X)^{-1} \t \left [ I(s(X))\right] I(t(X))^{-1}, \forall X \in A_1$,
\end{enumerate}  
is an $f$-homotopy. 
\end{Proposition}
\begin{Proposition}
(We continue to use the notation of the previous definition.)  Let $I$ be a
2-fold $f$-homotopy. Let $H$ be a homotopy with source $g$ and target $ f$. Then: 
\[(H\t I)(a)=H_1(a) \t I(a), \forall a \in A_0,\]
is  a 2-fold $g$-homotopy. Moreover, this assignment defines a left groupoid
action of $\CRS_1(\A,\G)$ on $\CRS_2(\A,\G)$ by automorphisms, and together with the boundary map $\d\colon   \CRS_2(\A,\G) \to \CRS_1(\A,\G)$,  already described, defines a crossed module $\CRS(\A,\G)$.  
\end{Proposition}

This is shown in \cite{BH4} in the general case when $\G$ is
a crossed complex. 
\subsection{Yetter's Invariant of  Manifolds}

\subsubsection{$\G$-colouring and Yetter's Invariant}\label{YI}
Let $\G=(G,E,\d,\t)$ be a finite reduced crossed module. Let $M$ be a triangulated piecewise linear
manifold or, more
generally,  a simplicial complex. We suppose that we have a total order on the set of vertices of
$M$.  Each simplex $K$ of $M$ has a representation in the form $K = \{ a,b,c,
 \ldots\}$, where $a<b<c< \ldots$ are vertices of $M$.  (Later, in section
 1.3.1, we will convert the ordered simplicial complex to a \emph{simplicial
 set} by using the total order.  The original simplices are then exactly the
 non-degenerate simplices of the result.)

\begin{figure}
\centerline{\relabelbox 
\epsfysize 3cm
\epsfbox{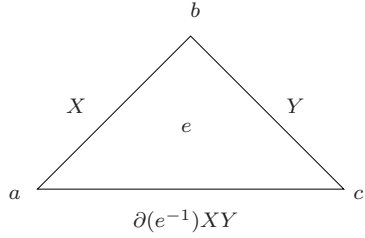}
\relabel {X}{$\s{X}$}
\relabel {Y}{$\s{Y}$}
\relabel {Z}{$\s{\d(e^{-1})XY}$}
\relabel {e}{$\s{e}$}
\relabel {1}{$\s{a}$} 
\relabel {2}{$\s{b}$}
\relabel {3}{$\s{c}$}
\endrelabelbox }
\caption{\label{simplex}A $\G$-colouring of an ordered $2$-simplex. Here $X, Y \in G,e \in E$ and $a<b<c$ are vertices of $M$.  To this $\G$-colouring we associate
  the  morphism  $XY \ra{e} \d(e)^{-1}XY$, in the category $\C(\G)$. }
\end{figure}

\begin{figure}
\centerline{\relabelbox 
\epsfysize 9cm
\epsfbox{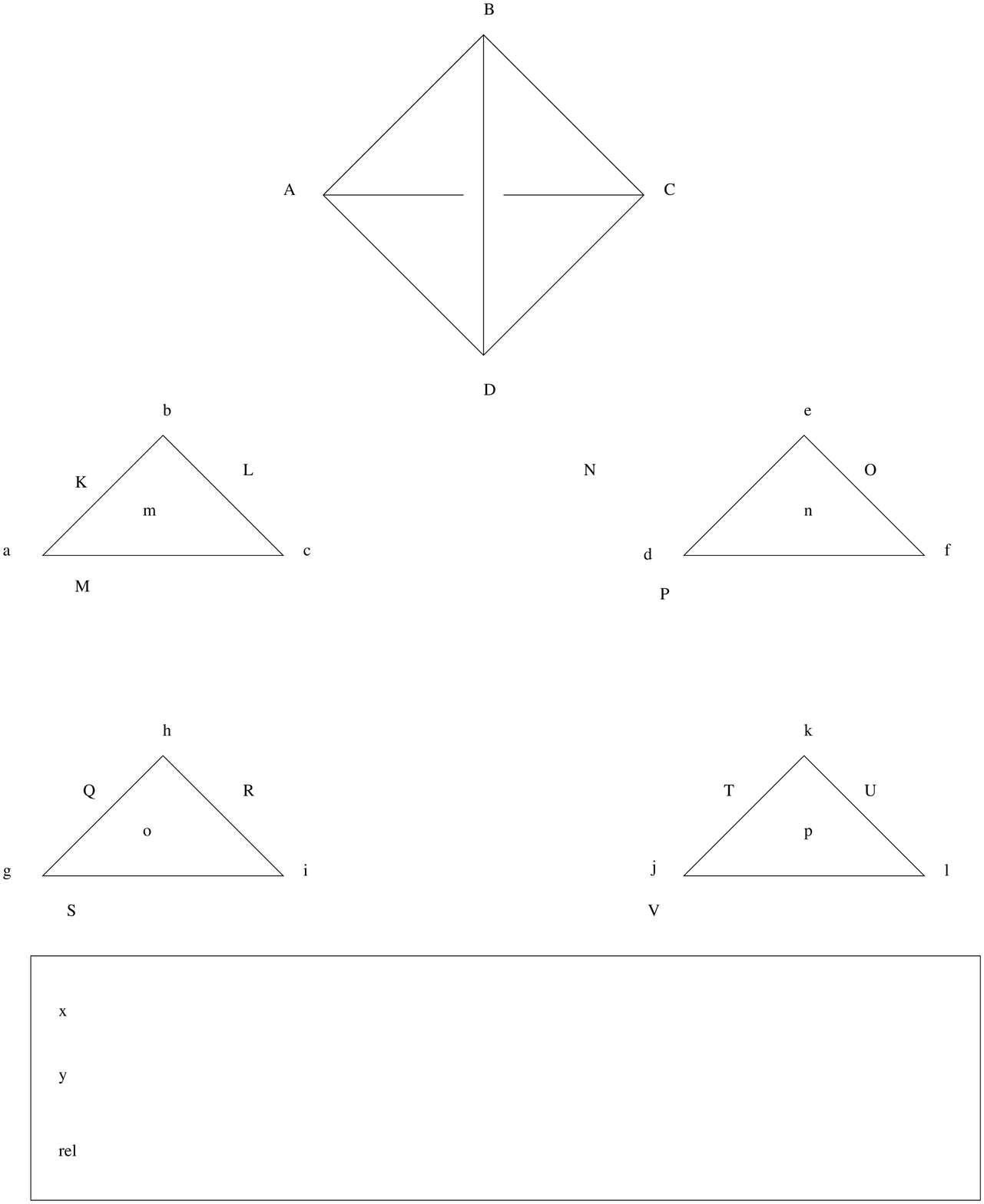}
\relabel {A}{$\s{a}$}
\relabel {B}{$\s{b}$}
\relabel {C}{$\s{c}$}
\relabel {D}{$\s{d}$}
\relabel {a}{$\s{a}$}
\relabel {b}{$\s{b}$}
\relabel {c}{$\s{c}$}
\relabel {d}{$\s{a}$}
\relabel {e}{$\s{c}$}
\relabel {f}{$\s{d}$}
\relabel {g}{$\s{b}$}
\relabel {h}{$\s{c}$}
\relabel {i}{$\s{d}$}
\relabel {j}{$\s{a}$}
\relabel {k}{$\s{b}$}
\relabel {l}{$\s{d}$}
\relabel {K}{$\s{X}$}
\relabel {L}{$\s{Y}$}
\relabel {M}{$\s{\d(e)^{-1}XY}$}
\relabel {N}{$\s{\d(e)^{-1}XY}$}
\relabel {O}{$\s{Z}$}
\relabel {P}{$\s{\d(f)^{-1}\d(e)^{-1}XYZ}$}
\relabel {Q}{$\s{Y}$}
\relabel {R}{$\s{Z}$}
\relabel {S}{$\s{\d(g)^{-1}YZ}$}
\relabel {T}{$\s{X}$}
\relabel {U}{$\s{\d(g)^{-1}YZ}$}
\relabel {V}{$\s{\d(h)^{-1}X\d(g)^{-1}YZ}$}
\relabel {x}{$\s{\psi=\{XYZ \ra{e} \d(e)^{-1}XYZ \ra{f} \d(f)^{-1} \d(e)^{-1} XYZ\}}$}
  \relabel {y}{$\s{\phi=\{XYZ \ra{X \t g} X\d(g)^{-1} YZ \ra{h} \d(h)^{-1} X
\d(g)^{-1}YZ\}}$}
\relabel {rel}{$\s{\phi=\psi \textrm{ is equivalent to } ef=(X \t g) h}$}
\relabel {m}{$\s{e}$}
\relabel {n}{$\s{f}$}
\relabel {o}{$\s{g}$}
\relabel {p}{$\s{h}$}
\endrelabelbox }
\caption{\label{tetrahedron}A $\G$-colouring of an ordered  tetrahedron. Here $a<b<c<d$.  Note that the
  compatibility condition $ef=(X \t g) h$ ensures that $\d(f)^{-1} \d(e)^{-1}
  XYZ= \d(h)^{-1} X \d(g)^{-1}YZ$. }
\end{figure}

A $\G$-colouring of $M$, also referred to as a ``Formal $\G$-Map'' in \cite{P3,PT}, is given by an assignment of an element of $G$ to
each edge of $M$ and an element of $E$ to each  triangle of $M$, satisfying the
compatibility condition shown in figure \ref{simplex}. A $\G$-colouring of a
 triangle naturally determines a morphism in $\C(\G)$, the monoidal category constructed from $\G$; see figure \ref{simplex}. There is still an
extra condition that $\G$-colourings must satisfy: Given a coloured
  3-simplex, one can
construct two morphisms $\f$ and $\p$ in
$\C(\G)$, see figure
\ref{tetrahedron}. Note that we need to use the monoidal structure of $\C(\G)$
to define $\f$ and $\p$. These morphisms are
required to be equal.

Let us be specific. Given a simplex $K$ of $M$, denote the colouring of it by
$\bc(K)$. Then the conditions that a $\G$-colouring  $\bc$ must satisfy  are:
\begin{equation}\label{colour2}
\d(\bc(abc))^{-1}\bc(ab) \bc(bc)=\bc(ac),
\end{equation}
whenever $(abc)$ (where $a<b<c$) is a $2$-simplex of $M$, and a cocycle condition:
\begin{equation}\label{colour3}
\bc(abc)\bc(acd)=\big (\bc(ab) \t \bc(bcd)\big )\bc(abd),
\end{equation}
whenever $(abcd)$ (where $a<b<c<d$) is a $3$-simplex of $M$.

\begin{Theorem} {\bf (D.N. Yetter)}
Let $M$ be a compact piecewise linear manifold. Consider a triangulation of $M$ with $n_0$ vertices
and $n_1$ edges. Choose a total order on the set of all vertices of $M$. Let
$\G=(G,E,\d,\t)$ be  a finite reduced crossed module. The quantity:
\begin{equation}\label{YetterInv}
I_\G(M)=\frac{\#E^{n_0}}{\#G^{n_0}\#E^{n_1}}{\# \left \{\G\textrm{-colourings of }
    M  \right \}}
\end{equation}
is a {piecewise linear homeomorphism} invariant of $M$, thus, in particular, it is triangulation independent, and does not depend on the total order chosen on the set of vertices of $M$.
\end{Theorem}
This result is due to D.N. Yetter.  Proofs are given   in \cite{Y2,P1,P2}. Similar
constructions  were considered  in \cite{FM1,FM2,FM3}, and also \cite{P3,PT}.
We will refer to the invariant $I_\G$, where $\G$ is a finite crossed module, as ``Yetter's Invariant''.

D.N. Yetter's proof of Theorem  \ref{main1} relied on Alexander Moves on
simplicial complexes, thereby  staying in the piecewise linear category. In particular it does not make clear whether $I_\G$ is a homotopy invariant.  On the other hand Yetter's argument extends $I_\G$ to being a TQFT.

In this article 
we will give a geometric interpretation of Yetter's Invariant $I_\G$, from  which will follow  a new
proof of its existence.  The expected fact that $I_\G$ is a homotopy
invariant is also implied by it; see \cite{P1,P2} for
alternative points of view, with a  broad intersection with this one.  As a
consequence, we  also have that
$I_\G(M)$ is well defined, as a homotopy invariant, for  any simplicial
complex $M$.

 In the closed piecewise linear  manifold
case, we will also consider an additional twisting of $I_\G(M)$ by cohomology
classes of $\G$, similar to the Dijkgraaf-Witten Invariant of 3-dimensional
manifolds, see \cite{DW}, thereby giving a solution to one of the open
problems posed by D.N. Yetter in \cite{Y2}.   Together with the homotopy invariance and geometric interpretation of Yetter's Invariant, this will be one of the main aims of this article.

\subsection{Classifying Spaces of Crossed Modules}

\subsubsection{The Homotopy Addition Lemma}\label{homadd}

Let $M$ be a piecewise linear manifold with a triangulation $T$, or, more
generally,  a simplicial complex. Suppose that we have a total
order on the set of all vertices of $T$.  Simplicial complexes like this will
be called ordered.

 Recall that we can define a simplicial set $D_T$ out of any
ordered  simplicial complex $T$. The set of $n$-simplices of $D_T$ is
given by all non decreasing sequences $(a_0,a_1,...,a_n)$ of vertices of $T$, such
that $\{a_0,a_1,...,a_n\}$ is a (non-degenerate) simplex of $T$. The face and degeneracy
maps are defined as:
\[
\d_i (a_0,...,a_n)=(a_0,...,a_{i-1},a_{i+1},...,a_n),i=0,...,n,
\]
and
\[
s_i (a_0,...,a_n)=(a_0,...,a_{i},a_{i},...,a_n),i=0,...,n,
\]
respectively. Notice that if $T$ is a piecewise linear triangulation of a piecewise linear manifold $M$ then the geometric realisation $|D_T|$ of $D_T$ is
(piecewise linear) {homeomorphic} 
to $M$.

 There exist several references on simplicial sets. We refer for example, \cite{C,FP,GZ,M,PK,We}.

 Working with simplicial sets, let $\Delta(n)$ be the  $n$-simplex, then its geometric realisation $|\Delta(n)|$ is the usual  standard geometric $n$-simplex.
 This simplicial set is
obtained from the simplicial complex $(012...n)$, the $n$-simplex, by applying the
construction above.   
 Each non-degenerate $k$-face $s=(a_0a_1a_2...a_k)$ of $\D(n)$ (thus $a_0<a_1<...<a_k$) 
 determines, in the obvious way, an element $h(s) \in
 \pi_k(|\D(n)|^k,|\D(n)|^{k-1},a_0)$. Recall that if $M$ is a CW-complex then
 $M^k$ denotes the $k$-skeleton of it.
 
By the General van Kampen Theorem, stated in \ref{general}, the crossed complex $\Pi(|\Delta(n)|)$ is
free on the set of non-degenerate faces of $\Delta(n)$; see
\cite{B1,B4,W4,BH5}. The definition of free crossed complexes appears for
example in \cite{B4}.  In particular, if $\A=(A_n)$ is a crossed complex, a morphism
$\Pi(|\Delta(n)|) \to \A$ is specified, uniquely,  by an assignment of an
element of $A_k$ to each non-degenerate $k$-face of $\Delta(n)$, with the obvious
compatibility relations with the boundary maps. To determine these relations,
we need to use the Homotopy Addition Lemma. This lemma explicitly describes
the boundary maps of the crossed complex $\Pi(|\Delta(n)|)$; see \cite{GW}, page
175. Let us be explicit in the low dimensional case. We follow the
conventions of \cite{BH5,B1}.

\begin{Lemma}{\bf (Homotopy Addition Lemma for $n=2,3$)}\label{HomAddLemma}
We have:
\begin{equation}
\d(h(0123))=\big(h(01)\t h(123)\big) h(013)h(023)^{-1}h(012)^{-1},
\end{equation}
\begin{equation}
\d(h(012))=h(12)h(23)h(13)^{-1}.
\end{equation}
\end{Lemma}
These equations are analogues  of equations (\ref{colour2}) and
(\ref{colour3}).  
By the  General van Kampen Theorem  there follows:
\begin{Proposition}\label{Refer1}
Let $M$ be a triangulated piecewise linear manifold with a total order on its set of vertices. Let also $\G$ be a reduced crossed
module. There exists a one-to-one correspondence between $\G$-colourings of
$M$ and
maps $\Pi(M) \to \G$, where $M$ has the natural CW-decomposition given by its
triangulation. 
\end{Proposition}
This result appears in \cite{P1} (Proposition 2.1). Note that any map $f\colon   \Pi(M) \to \G$ yields
a $\G$-colouring $\bc^f$ of $M$, where $\bc^f(K)=f(h(K))$. Here $K$ is a non-degenerate simplex of $M$. To prove this we need
to use the Homotopy Addition Lemma.  

Proposition \ref{Refer1}   still holds for general simplicial complexes, and has obvious analogues if we merely specify a CW-decomposition of $M$ rather than a triangulation.

More discussion of   $\G$-colourings can be found  in \cite{P1,P2,P3,PT}.

\subsubsection{The Definition of $B(\G)$ and Some Properties}\label{Refer2}
We follow \cite{BH5,B1,B4}. Another introduction to this subject appears in \cite{P3}.
Let $\A$ be a crossed complex.  The nerve  $N(\A)$  of $\A$ is, by
definition, the simplicial set  whose set of $n$-simplices is given by all
crossed complex morphisms $\Pi(|\D(n)|) \to \A$, with the obvious  face and degeneracy maps.   It is a Kan simplicial set; see \cite{BH5}. The classifying space $B(\A)$ of $\A$ is the geometric
realisation  $|N(\A)|$ of the nerve of $\A$.  This construction appeared first
in \cite{Bl}.

Let us unpack this definition in the reduced crossed module case; see
\cite{B1}, $3.1$. Let $\G=(G,E,\d,\t)$ be a reduced crossed module. The CW-complex $B(\G)$
has a unique $0$-cell $*$, since $N(\G)$  only has one $0$-simplex, because $\G$ is reduced.  For any $n\in \N$, the
set of $n$-simplices of $N(\G)$ is, by Proposition \ref{Refer1}, in bijective correspondence with  the set of all
$\G$-colourings of $|\D(n)|$, and has the obvious faces and degeneracies motivated by this identification.

Consequently, there exists a $1$-simplex $s(X)$ of $B(\G)$ for
any $X \in G$. The 2-simplices of $M$ are given by triples $s(X,Y,e)$, where
$X,Y \in G$ and $e \in E$. Here $X=\bc(01),Y=\bc(12),\d(e)^{-1}XY=\bc(02)$ and
$e=\bc(012)$, corresponding to the $\G$-colouring of $|\D(2)|$ shown in figure
\ref{simplex}.  Similarly, the 3-simplex of $N(\G)$ determined by the $\G$-colouring of  $|\D(3)|$ shown in figure
\ref{tetrahedron}  will be denoted by $s(X,Y,Z,e,f,g,h)$, where the relation
$ef=(X \t g) h$ must hold.  We will describe  the set of $\G$-colourings of the
4-simplex $|\D(4)|$ in the following section.

The following result is due to R. Brown and P.J. Higgins; see \cite{BH5},
Theorem $2.4$, together with its proof.

Let $C$ be a simplicial set. In particular the geometric realisation $|C|$ of
$C$ is a CW-complex.
\begin{Theorem}\label{Refer4}
 Let $\G$ be a crossed module. There exists a one-to-one
correspondence $F$ between crossed complex morphisms $\Pi(|C|)\to \G$ and
simplicial maps $C \to N(\G)$. 
\end{Theorem}
\begin{Remark} \label{REFER}
 In fact,
 let $g\colon    \Pi(|C|) \to \G$ be a
morphism of crossed complexes. Consider the projection map: \[p\colon    \bigsqcup_{n \in
  \N}C_n \times |\D(n)| \to |C|,\]
according to the usual construction of the geometric realisation of simplicial
sets, see for example \cite{Mi1,M,FP,GZ}. Note that $p$ is a cellular map,
therefore it induces a morphism $p_*$ of crossed complexes. 
  Let  $c\in C_n$ be  an $n$-simplex of $C$, then $F(g)(c)\colon    \Pi(|\D(n)|) \to \G$ is
the composition: 
\[\Pi(|\D(n)|)\ra{\cong} \Pi(c\times |\D(n)|)\ra{p_*} \Pi(|C|) \ra{g} \G.\]
\end{Remark}
These results are a  consequence of the version of  the General van
 Kampen Theorem stated in \ref{general}; see \cite{BH5} for details.
Remark \ref{REFER} shows  how to go from a crossed complex
morphism to a simplicial morphism.  The reverse direction needs a
comparison of $\Pi (B(\mathcal{G}))$ and $\mathcal{G}$ and is more
complicated.   However,  we can prove directly that  $F$ is a bijection by using the General van Kampen Theorem.

\begin{Remark}\label{ColSimp}
In   particular (from Proposition \ref{Refer1}), it follows that if $M$ is a piecewise linear  manifold with an
ordered triangulation $T$, and
$C_T$ is the simplicial set associated with $T$ (see the beginning of \ref{homadd}), then there
exists a one-to-one correspondence between $\G$-colourings of $M$ and
simplicial maps $C_T\to N(\G)$. This result is due to the second author; see \cite{P1,P2}. 
\end{Remark}
The statement of Theorem \ref{Refer4} can be substantially expanded. We follow \cite{BH5}. The assignment
$\G \mapsto N(\G)$ if $\G$ is a crossed module is functorial. Similarly, the
assignment $C \mapsto \Pi(|C|)$, where $C$ is a simplicial
set is also functorial.  The second functor is left adjoint to the first one:
the Nerve Functor,  see \cite{BH5}, Theorem 2.4. In fact, the adjunction $F$ of
Theorem \ref{Refer4} extends to  a simplicial map (which we also call $F$): 
\[N(\CRS(\Pi(|C|),\G)) \ra{F} \SIMP(C,N(\G)),\]
which, moreover, is  a homotopy equivalence; see \cite{BH5}, proof of Corollary $3.5$ and of  Theorem $A$. 
Here if $C$ and $D$ are simplicial sets, then $\SIMP(C,D)$ is the simplicial set
whose set of $n$-simplices is given by all simplicial maps $C\times \D(n) \to D$; see
\cite{M,GZ}, for example. Note that the set of vertices of
$N(\CRS(\Pi(|C|),\G))$ is the set $\CRS_0(\Pi(|C|),\G)$ of all crossed complex morphisms $\Pi(|C|)
\to \G$. Similarly, there exists a bijective correspondence between vertices of
$\SIMP(C,N(\G))$ and simplicial maps $C \to N(\G)$.

The existence of the simplicial  homotopy equivalence $F$ is a non  trivial fact, due to
 R. Brown and P.J. Higgins. The  proof which appears in \cite{BH5} makes great  use  of the monoidal closed structure of the category of crossed complexes constructed in \cite{BH4}.  A version of the  Eilenberg-Zilber Theorem for crossed complexes, representing $\Pi(|C \times D|)$ as a strong deformation retract of $\Pi(|C|)\otimes \Pi(|D|)$, if $C$ and $D$ are simplicial sets, is also required.  An approach to this is  sketched in \cite{BH5}. A direct proof appears in  \cite{To1,To2}.

Let $C$ and $D$ be simplicial sets. It is well known that if $D$ is Kan then there exists a weak homotopy equivalence: 
\[j\colon   |\SIMP(C,D)| \to \TOP(|C|,|D|), \]
where the later is given the $k$-ification of the Compact-Open Topology on the
set of all continuous maps $|C| \to |D|$; see \cite{B2}, for example.  This weak homotopy equivalence is the composition:
\[|\SIMP(C,D)| \ra{j_1}|\SIMP(C,S|D|)| \ra{j_2} |S(\TOP(|C|,|D|))| \ra{j_3} \TOP(|C|,|D|).\]
Here if $M$ is a topological space, then $S(M)$ is the singular complex of
$M$. Since  $D$ is a Kan simplicial set, it follows that $D$ is a strong deformation 
retract of $S(|D|)$,  a fact usually known as Milnor's Theorem, see \cite{FP},
$4.5$, and  \cite{Mi1}.  The map $j_1$ is the geometric realisation of the inclusion
$\SIMP(C,D)\to \SIMP(C,S|D|)$, and is therefore a homotopy equivalence.

There exists a one-to-one correspondence between simplicial maps $C \times \D(n) \to
S(|D|)$ and continuous maps $|C| \times |\D(n)| \cong |C \times \D(n)| \to |D|$. This fact together
with the natural homeomorphism    $\TOP(|C| \times |\D(n)|,|D|) \to  \TOP(|\D(n)|,
\TOP(|C|,|D|))$ provides an  isomorphism $\SIMP(C,S(|D|)) \to
S(\TOP(|C|,|D|))$; see \cite{GZ}, page 132. The map $j_2$ is the geometric
realisation of it, which consequently is  a homeomorphism.

The map $j_3$ is the obvious map
$j_3\colon   |S(\TOP(|C|,|D|))| \to \TOP(|C|,|D|)$. It is well known that the natural
map $|S(M)| \to M$ is a weak homotopy equivalence, for any space $M$; see
\cite{M,GZ,FP}. 

 From this discussion, it is obvious that if $f\colon   C \to D$ is a
simplicial map then $j(|f|)=|f|$. Here $|f|\colon   |C| \to |D|$ is the geometric
realisation of $f\colon   C \to D$. 
Therefore:
 
\begin{Theorem}\label{weak}
Let $C$ be a simplicial set and $\G$ a reduced crossed module. There exists a weak
 homotopy equivalence: 
\[\n\colon   B(\CRS(\Pi(|C|),\G)) \to \TOP(|C|,B(\G)).\]
In fact, if we are given a map  $g\colon   \Pi(C) \to \G$, then $\n(g)$ is the
 geometric realisation of $F(g)\colon    C
 \to N(\G)$; see Theorem \ref{Refer4}.
\end{Theorem}
This is  contained in theorem $A$ of \cite{BH5}. Note that $\n=j \circ |F|$. The exact expression of $\n$ will be needed in the last chapter.

We  therefore have that $\n$ is  a homotopy equivalence when $C$ is finite, because then it follows
that $\TOP(|C|,B(\G))$ has the homotopy type of a CW-complex; see \cite{Mi2}.

\subsubsection{The Fundamental Groups of Classifying Spaces of Crossed Modules}

Let $\G=(G,E,\d,\t)$ be a crossed module, not necessarily reduced. Let $P$ be the object set of $G$ and $E$. There exists a bijective correspondence between 0-cells of $B(\G)$ and elements of $P$. 

\begin{Theorem}\label{fundamental}
There exists a one-to-one correspondence between connected components of
$B(\G)$ and connected components of the groupoid $G$. Given $a\in P$ then:
\begin{equation}
\pi_1(B(\G),a)=G_{a,a}/{\rm im} \{\d\colon   E_a \to G_{a,a}\}
\end{equation} 
and
\begin{equation}
\pi_2(B(\G),a)=\ker \{\d\colon   E_a \to G_{a,a}\}.
\end{equation} 
Here $G_{a,b}=\{X \in G\colon   s(X)=a,t(X)=b\}$ and $E_a=\{e \in E\colon   \b(e)=a\}$, where
$a,b \in P$. 
In fact the groupoid $\pi_1(B(\G),P)$ is the quotient of $G$ by the totally
disconnected subgroupoid $\d(E)$, which is normal in $G$. All the  other homotopy groups of $B(\G)$ are trivial. 
\end{Theorem}
 This is shown in \cite{BH5}, in the general case of crossed complexes. In fact, if $\A=(A_n)$ is a {reduced} crossed complex then $\pi_n(B(\A))=\ker ( \d_n) /\mathrm{im} (\d_{n+1})$, {and analogously in the non-reduced case, keeping track of base points.}   These results are a   consequence of the fact that the nerve $N(\A)$ of $\A$ is a Kan simplicial set if $\A$ is a crossed complex.

 An immediate corollary of theorems \ref{weak} and \ref{fundamental} is that
if $C$ is a simplicial set then there exists a one-to-one correspondence
between homotopy classes of crossed complex morphisms $\Pi(|C|) \to \G$ and
homotopy classes of maps $|C| \to B(\G)$. More precisely, $f,g\colon    \Pi(|C|) \to
\G$  are homotopic if and only if $\n(f),\n(g)\colon    |C| \to B(\G)$ are
homotopic. R. Brown and P.J. Higgins proved this for any CW-complex, a statement contained in Theorem $A$ of \cite{BH5}.

\subsubsection{n-Types, {in General, and 2-Types, in Particular}}

Let $\A$ be a crossed complex. The assignment $\A \mapsto N(\A)$, {where $N(\A)$ is the nerve of $\A$}, is functorial. Composing with the geometric realisation functor from the category of simplicial sets to the category ${Top}$ of topological spaces, proves that the geometric realisation map  $\A \mapsto B(\A)$ can be extended to a functor ${Xcomp} \to {Top}.$ Here  ${Xcomp}$ is the category of crossed complexes.  

Recall that a continuous map $f : X\to Y$ is an \emph{$n$-equivalence} if
it induces a bijection $\pi_0(f) $, on connected components and if, for
all choices of base point $x_0\in X$ and $1\leq i \leq n$, $\pi_i(f) :
\pi_i(X,x_0) \to \pi_i(Y,f(x_0))$ is an isomorphism.  (Intuitively an
$n$-equivalence is a truncated weak homotopy equivalence as it  ignores
all information in the $\pi_k(X)$ for $k> n$.) Two spaces have the same
$n$-type if there is a zig-zag of $n$-equivalences joining them, more
exactly one can formally invert $n$-equivalences to get a category
$Ho_n(Top)$ and then $X$ and $Y$ have the same $n$-type if they are
isomorphic in $Ho_n(Top)$.  Any CW-space $X$ has a $n$-equivalence to a
space $Y_n$ (the $n$-type of $X$), which satisfies $\pi_k(Y_n)$ is trivial for all $k>n$.  This
$Y_n$ is formed by adding high dimensional cells to $X$ to kill off the
elements of the various $\pi_k(X)$.  Because of this, it is usual,
loosely, to say that a CW-complex $M$ \emph{is an $n$-type} if $\pi_i(M)=0$
for all $i>n$.  We will adopt this usage. Note that the $n$-type of a CW-complex is well defined up to homotopy.

If $\A=(A_m)$ is a crossed complex such that $A_m$ is
trivial  for $m>n$ then it follows that  the classifying space $B(\A)$ of $\A$
is an $n$-type. The reader is advised that the assignment $\A \mapsto B(\A)$ where $\A$ is a crossed complex does not generate all $n$-types.  

Working with reduced crossed modules, it is well known that any  connected $2$-type can be represented as the classifying space of a reduced crossed module,  a result which is essentially due to S. Mac Lane and J.H.C Whitehead; see \cite{MLW}. See
for example \cite{BH5} for an alternative proof. In fact it is possible to prove that the {homotopy} category of {connected} 2-types is equivalent to $Ho(\mathrm{X}-mod)$, where $\mathrm{X}-mod$ is the category of reduced crossed modules; see \cite{BA1,BH5,L}. Here,  a weak equivalence of reduced crossed modules is by definition a map $f{\colon}\G \to {\cal H}$ whose geometric realisation induces isomorphisms $\pi_1(B(\G),*)\to  \pi_1(B({\cal H}),*)$ and  $\pi_2(B(\G),*)\to  \pi_2(B({\cal H}),*)$. By Theorem \ref{fundamental}, this can be stated in algebraic terms. 
 For the description of the associated model structure on the category of crossed modules, we refer the reader to \cite{BG}. 

 In particular, it follows that  the reduced crossed modules $\G$ and ${\cal H}$ are weak equivalent if and only if $B(\G)$ and $\B({\cal H})$ have the same $2$-type. Since these CW-complexes  are 2-types {themselves},  there follows that:
\begin{Theorem}\label{A}
Let $\G$ and ${\cal H}$ be reduced crossed modules. Then $\G$ and ${\cal H}$ are weak homotopic if and only if $B(\G)$ and $B({\cal H})$ are homotopic.
\end{Theorem}

Let $\G=(G,E,\d,\t)$ be a crossed module. There exists a natural action of $\mathrm {coker} (\d)$ on $\ker( \d)$.  The crossed module $\G$ determines a cohomology class $\w \in H^3(\ker (\d),\mathrm{coker}(\d))$, called the $k$-invariant of $\G$, which coincides with the $k$-invariant of $B(\G)$ (also called first Postnikov invariant); see \cite{B1}. By a theorem of Mac Lane and Whitehead (see \cite{MLW}), two spaces have the same $2$-type if and only if they have the same fundamental group, the same second homotopy group, seen as a $\pi_1$-module, and the same {first} Postnikov invariant; see \cite{BA1,BA2}. 

A triple $(A,B,k^3 \in H^3(A,B))$, where $A$ is a group acting on the abelian group $B$ on the left by automorphisms  is usually called an algebraic 2-type. The discussion above implies:

\begin{Corollary}\label{B}
Let ${\G}$ and ${\cal H}$ be crossed modules. Then $B(\G)$ is homotopic to $B({\cal H})$ if and only if the algebraic 2-types determined by $\G$ and ${\cal H}$  are isomorphic. 
\end{Corollary} 

\subsection{Applications to Yetter's Invariant of Manifolds}
\subsubsection{An Interpretation of Yetter's Invariant}
Let $\G$ be a finite  reduced crossed module. Recall that $I_\G$ stands for  Yetter's Invariant; see \ref{YI}. We freely use the results and the notation of \ref{HomMaps}.
We now prove the first main result of this article.

\begin{Theorem}\label{main1}
Let $M$ be a {compact} triangulated piecewise linear manifold or, more generally,  a {finite} simplicial complex.  As usual, we suppose that $M$ is provided with a total order on its set of vertices. Let $\G=(G,E,\d,\t)$ be a reduced crossed module. We have:
\begin{equation}\label{D}
I_\G(M)=\sum_{g \in [M,B(\G)]} \frac{\#\pi_2(\TOP(M,B(\G)),g)}{\#\pi_1(\TOP(M,B(\G)),g)}.
\end{equation}
Here $[M,B(\G)]=\pi_0(\TOP(M,B(\G))$ denotes the set of homotopy classes of maps $M \to B(\G)$. In
particular, $I_\G(M)$ does not depend on the ordered triangulation of $M$ chosen. Moreover, Yetter's Invariant is therefore a homotopy invariant of polyhedra, i.e, triangulable topological spaces.
\end{Theorem}
{Note that equation (\ref{D}) extends Theorem 12 of \cite{Y2} to general crossed modules}

\begin{Proof}
The two previous theorems will be used. Recall that any simplicial complex
naturally defines a simplicial set; see \ref{homadd}. We freely use the
notation and the results of \ref{HomMaps}. 
We  abbreviate $\CRS_2(\Pi(M),\G)(f)$ as $\CRS_2(f)$, and analogously
for $\CRS_1(\Pi(M),\G)(f,g)$ and  $\CRS_1(\Pi(M),\G)(f)$.

For an $f \in \CRS_0(\Pi(M),\G)$,  denote  $[f]$ as being the set of elements of
$\CRS_0(\Pi(M),\G)$ which belong to  the same connected component of $f$ in
the groupoid $\CRS_1(\Pi(M),\G)$. Recall  also that if $f,g\colon   \Pi(M) \to \G$ then
$\n(f),\n(g)\colon   M \to B(\G)$  are homotopic if, and only if, $f$ and $g$ belong
to the same connected component in the groupoid $\CRS_1(\Pi(M),\G)$.  

  We have:
\begin{align*}
\sum_{g \in [M,B_\G]} &\frac{\#\pi_2(\TOP(M,B(\G)),g)}{\#\pi_1(\TOP(M,B(\G)),g)}
\\&=\sum_{f \in \CRS_0(\Pi(M),\G)} \frac{\#\pi_2(\TOP(M,B(\G)),\n(f))}{\#\pi_1(\TOP(M,B(\G)),\n(f)) \#[f]}
\\&=\sum_{f \in \CRS_0(\Pi(M),\G)} \frac{\#\pi_2(B(\CRS(\Pi(M),\G)),f)}
{\#\pi_1  (B(\CRS(\Pi(M),\G)),f) \#[f]}
\\ &=\sum_{f \in \CRS_0(\Pi(M),\G)} \frac{\# \d(\CRS_2(f))\# \ker \{\d\colon   \CRS_2(f) \to \CRS_1(f,f)\}  }{\#\CRS_1(f,f) \#[f]}
\\  &=\sum_{f \in \CRS_0(\Pi(M),\G)} \frac{\# \CRS_2(f)}{\#\CRS_1(f,f) \#[f]}.
\\  &=\sum_{f \in \CRS_0(\Pi(M),\G)} \frac{\# \CRS_2(f)}{\#\CRS_1(f)}.\\
\end{align*}
The last equality follows from the fact that $\CRS_1(\Pi(M),\G)$ is a
groupoid. Recall that $\CRS_1(f)$ denotes the set of homotopies with source
$f\colon\Pi(M) \to \G$.
The result follows from the following  lemmas. 
\begin{Lemma}
Let $\f\colon   H \to G$ be a groupoid morphism, with $G$ a group. Let also $E$ be a
group on which $G$ acts on the left. Suppose that $H$ is a free groupoid. Any $\f$-derivation $d\colon   H \to E$ can be specified, uniquely, by its value on the
set of free generators of $H$. In other words, any map from the set of
generators of $H$ into $E$ uniquely extends to a derivation $d\colon H \to E$. 
\end{Lemma}
\begin{Proof}
This is shown in \cite{W4} (Lemma $3$) for the case in which $H$ is a group, and is 
easy to prove. The general case can be proved directly in the same
way. Alternatively, we can reduce it to the group case by  using the free
group $\hat{H}$ on the set of free generators of $H$. In particular,  any
groupoid morphism $\f\colon   H \to G$ factors uniquely as $\f=\p \circ p$, where $p$ is the obvious map $p\colon   H\to
\hat{H}$, bijective on generators, and $\psi\colon \hat{H} \to G$ is a group morphism. Note that if $d\colon   \hat {H}  \to E$ is a
$\p$-derivation, then  $d \circ p\colon   H \to E$ is a $\f$-derivation.   
\end{Proof}

\begin{Lemma}
 There exist exactly $\#E^{n_1}\#G^{n_0}$ $f$-homotopies if $f\colon   \Pi(M) \to \G$ is a crossed complex morphism. Here $n_0$ and $n_1$ are the number of vertices and edges  of $M$.
\end{Lemma}
\begin{Proof}
Follows from the previous lemma and the fact that $\pi_1(M^1,M^0)$ is the free groupoid on the set of edges of $M$, with one object for each vertex of $M$. The fact that the crossed module $\G=(G,E,\d,\t)$ is reduced is being used.
\end{Proof}
 \begin{Lemma}
There exist exactly $\#E^{n_0}$ 2-fold $f$-homotopies if $f\colon   \Pi(M) \to \G$ is a crossed complex morphism.
\end{Lemma}
\begin{Proof}This is a consequence of   the fact that the crossed module $\G$ is, by assumption,  reduced. 
\end{Proof}

Using the last two lemmas it follows that:
\begin{align*}
\sum_{f \in [M,B(\G)]}& \frac{\#\pi_2(\TOP(M,B(\G)),f)}{\#\pi_1(\TOP(M,B(\G)),f)}
\\ &=\sum_{f \in \CRS_0(\Pi(M),\G)} \frac{\# \CRS_2(f)}{\#\CRS_1(f)}
\\ &=\sum_{f \in \CRS_0(\Pi(M),\G)} \frac{\# E^{n_0}}{\# G^{n_0} \# E^{n_1}}
\\ &=\frac{\# E^{n_0}}{\# G^{n_0} \# E^{n_1}} \# \{\G\textrm{-colourings of } M\}.
\end{align*}
The last equation follows from Proposition \ref{Refer1}. This finishes the proof of Theorem \ref{main1}.
\end{Proof}
 
\begin{Corollary}
The invariant $I_\G$ where $\G$ is a finite crossed module depends only on the homotopy type of $B(\G)$. Therefore Yetter's invariant $I_\G$ is a function of the weak homotopy type of the crossed module $\G$, only, or alternatively of the algebraic 2-type determined by $\G$.
\end{Corollary}
\begin{Proof}
Follows from Theorem \ref{A} and Corollary \ref{B}.
\end{Proof}

In a future work, we will investigate the behaviour of $I_\G$ under arbitrary crossed module maps $\G \to {\cal H}$, in order to  generalise this result.

Yetter's Invariant was extended to CW-complexes in \cite{FM2}, where  an alternative proof of its existence is given,  using a totally different argument to the one just shown.  
Also in the cellular category, an extension of Yetter's Invariant to crossed
complexes appeared in \cite{FM3}. In \cite{FM1,FM2}, algorithms are given for
the calculation of Yetter's Invariant of complements of knotted embedded surfaces in $S^4$, producing non-trivial invariants of 2knots.
Different interpretations of Yetter's Invariant appear in \cite{P1,P2}.

{Yetter's Invariant has a very simple expression in the case of spaces homotopic to 2-dimensional complexes. In fact:}
\begin{Proposition} 
Let $\G=(G,E,\d,\t)$ be a finite reduced crossed module. Suppose that $M$ is a finite simplicial complex which is homotopic to a  finite 2-dimensional simplicial complex. Then:
\[I_\G(M)=\frac{1} {\# \mathrm{coker} (\d)}\big (\#\mathrm{ker}(\partial) \big)^{\chi(M)}  \# \mathrm{Hom}\big (\pi_1(M), \mathrm{coker} (\partial)\big ),\]
wbere $\chi(M)$ denotes the Euler characteristic of $M$.
\end{Proposition}
\begin{Proof} 
{It suffices to prove the result for  2-dimensional simplicial complexes, since Yetter's Invariant is an invariant of homotopy types. }

{
Let $\cal H$ be the crossed module with base group $ \mathrm{coker} (\partial)$, whose  fibre group is the trivial group.  There exist a natural map $p$ from the set of $\G$-colourings of $M$ to the set of $\cal H$-colourings of $M$, obtained by sending the colour of an edge of $M$ to the image of it in $\mathrm{coker} (\partial)$ and the colour of each face to $1$. Since there are no relations motivated by tetrahedra, the map $p$ is surjective, and the inverse image at each point has cardinality $(\#\mathrm{im}(\d))^{n_1}(\# \ker( \d))^{n_2}$. Here $n_i$ is the number of non-degenerate $i$-simplices of $M$. Therefore:} {
\begin{align*}
I_\G(M)&=\frac{\# E^{n_0}} {\# G^{n_0} \#E^{n_1}}(\# \mathrm{im}(\d))^{n_1}(\# \ker( \d))^{n_2}\#\{{\cal H}\textrm{-colourings of } M\}\\
&=\frac{\# E^{n_0}} {\# G^{n_0}}(\# \ker( \d))^{n_2-n_1}\#\{{\cal H}\textrm{-colourings of } M\}\\
&=\frac{(\# \mathrm{im}(\d))^{n_0}} {\# G^{n_0}}(\# \ker( \d))^{n_2-n_1+n_0}\#\{{\cal H}\textrm{-colourings of } M\}\\
&=\frac{1} {(\# \mathrm{coker} (\d))^{n_0}}(\# \ker( \d))^{n_2-n_1+n_0}\#\{{\cal H}\textrm{-colourings of } M\}\\
&=\frac{1} {\# \mathrm{coker} (\d)}(\# \ker( \d))^{\chi(M)}\# \mathrm{Hom}\big (\pi_1(M), \mathrm{coker} (\partial)\big ).
\end{align*} }  
\end{Proof}
\subsubsection{An Extension of Yetter's Invariant to Crossed Complexes}\label{ABC}

We can extend Yetter's Invariant of manifolds to handle crossed complexes. Let us sketch how to do this extension. For details (in the pointed CW-complex category) we refer the reader to \cite{FM3}. 

As part of the information used to construct this extension, we need to know
the $n>3$ cases  of the Homotopy Addition Lemma, stated for $n=2$ and $n=3$ in
\ref{homadd}. We continue to follow the conventions of \cite{BH5}. 

\begin{Lemma} {\bf (Homotopy Addition Lemma for $n>3$)}
We have:
\[\d h(012...n)=h(01)\t h\big(\d_0 (012...n)\big)+\sum_{i=1}^n (-1)^i h\big (\d_i(012...n)\big).\]

 \end{Lemma}
Recall the notation introduced in \ref{homadd}. Note that this is an identity in the relative homotopy group $\pi_n(|\D(n)|^n,|\D(n)|^{n-1},0=*)$. 
\begin{Definition}
Let $M$ be a triangulated  piecewise linear manifold with a total order on its set of
vertices. Let $\A=(A_i)$ be a reduced crossed complex. An $\A$-colouring of $M$ is given by  an
assignment of an element of $A_i$ to each $i$-simplex of $M$, which should
satisfy the obvious relations motivated by the Homotopy Addition Lemma.
\end{Definition}
As before, by the General van Kampen Theorem, it follows that there is a one-to-one correspondence between $\A$-colourings of $M$ and crossed complex  maps $\Pi(M) \to \A$; compare with Proposition \ref{Refer1}. Using the same techniques that we used in the crossed module case we can prove:
\begin{Theorem}\label{Main3}
Let $M$ be a triangulated compact piecewise linear manifold (or more
generally a {finite} simplicial complex)  provided with a total order on its set of vertices, having  $n_i$
(non-degenerate) simplices of dimension $i$ for each $i \in \N$. Let also
$\A=(A_n)$ be a  reduced crossed complex such that $A_n$ is the trivial group except for a finite number of $n\in \N$ for which $A_n$ should be  finite. The quantity:
\[I_\A(M)=\#\{\A\textrm{-colourings of } M\}\prod_{i=1}^{\infty}\left (\prod_{j=0}^{\infty} (\# A_{i+j})^{n_j}\right )^{{(-1)}^{i}},\]
does not depend on the ordered triangulation of $M$ chosen and it is a homotopy invariant of
$M$, as a topological space. In fact $I_\A(M)$ can be interpreted as:
\[I_\A(M)=\sum_{f \in \pi_0(\TOP(M,B(\A)))} \prod_{i=1}^{\infty} \# \pi_i\big(\TOP(M,B(\A)), f\big)^{(-1)^i}.\]
Here $B(\A)$ is the classifying space of the crossed complex $\A$.
\end{Theorem}
The proof of this theorem is analogous to the proof of Theorem \ref{main1}. An important fact to note is  that if $\A$ and $M$ are  as in Theorem
\ref{main1}, then the function space $\TOP(M,B(\A))$ has only a finite number of connected
components, each with a finite number of non-trivial homotopy groups, all of
which are finite. This follows from  theorems \ref{weak} and   \ref{fundamental}, or more precisely from their extension to the general case of crossed complexes. A proof of Theorem
\ref{Main3} appears in \cite{FM3}, in the case when $M$ is a CW-complex with a
unique 0-cell.  Using the framework we exposed here in the crossed
module case, the proof is easily extended to the case when $M$ may have any finite number of $0$-cells; see the theorem below.

This construction should be compared with the  construction in \cite{P2} of a
manifold invariant from any finite $n$-type. It would be interesting to formulate  a similar homotopy interpretation of it. 

Theorems \ref{main1} and \ref{Main3} can be stated in a  more   general form in the CW-category. In fact we have:
\begin{Theorem}
Let $M$  be a finite CW-complex, having  $n_i$ cells  of index $i$ for each $i \in \N$. Let also
$\A=(A_n)$ be a  reduced crossed complex, such that $A_n$ is the trivial group, except for a finite number of $n\in \N$, for which $A_n$ should be  finite. The (necessarily finite) quantity:
\[I_\A(M)=\#\{ \rm{Hom}( \Pi(M), \A)\}\prod_{i=1}^{\infty}\left (\prod_{j=0}^{\infty} (\# A_{i+j})^{n_j}\right )^{{(-1)}^{i}},\] 
 is a homotopy invariant of
$M$, as a topological space. Here $\Pi(M)$ is the fundamental crossed complex of the skeletal filtration of $M$. Moreover  $I_\A(M)$ can be interpreted as:
\[I_\A(M)=\sum_{f \in \pi_0(\TOP(M,B(\A)))} \prod_{i=1}^{\infty} \# \pi_i\big(\TOP(M,B(\A)), f\big)^{(-1)^i},\]
where $B(\A)$ is the classifying space of the crossed complex $\A$. 
\end{Theorem}
This theorem is shown in \cite{FM3}. The proof is set in the based case. However the argument extends easily to the unbased case by the method shown in the crossed module {context}. 

A more combinatorial version of the previous theorem can always be set by using the General van Kampen Theorem, as long as we can apply an appropriate homotopy addition lemma, similar to the one stated in the simplicial context. This can be done for instance for complements of knotted embedded surfaces in $S^4$; see \cite {FM2}.
\section{Homological Twisting of Yetter's Invariant}
\subsection{Low Dimensional Cohomology of Crossed Modules}
The (co)homology of crossed modules is studied, for example, in \cite{E,Pa},
by considering the (co)homology of their classifying spaces. We do not wish to
develop that further. Rather, we will only consider  the 
(co)homology in low dimensions, focusing on a geometric and calculational  approach. Our final aim is to define a 
twisting of Yetter's Invariant by cohomology classes of finite reduced crossed
modules, thereby extending the Dijkgraaf-Witten Invariant of manifolds to crossed modules (equivalent to  categorical groups).

\begin{Definition}
Let $\G$ be a  reduced crossed module. The (co)homology of $\G$ is defined as being the (co)homology of the classifying space $B(\G)$ of $\G$, and similarly in the non-reduced case and for crossed complexes.
\end{Definition}

\subsubsection{Homology of Simplicial Sets}\label{homKan}
Let $D$ be a simplicial set. It is well known  that  the  homology groups  of
its geometric realisation $|D|$ can be  combinatorially defined from $D$
itself. Let us clarify this statement.

Consider the complex $C(D)=\{C_n(D), \d_n\}$ of simplicial chains of
$D$. Here $C_n(D)$ is the free abelian group on the set $D_n$ of $n$-simplices
of $D$. Furthermore $\d(c)=\sum_{i=0}^n (-1)^i \d_i(c)$, if $c \in D_n$; see
\cite{M}, $1.2$, or \cite{We}, $8.2$.  Note that the assignment $D \mapsto C(D)$,
where $D$ is a simplicial set is  functorial. The homology of a simplicial set
$D$ is defined as being  the homology of the chain complex $C(D)$.

The chain complex $C(D)$ has a
subcomplex $C^d(D)$, for which $C^d_n(D)$ is the free $\Z$-module on the set of
degenerate $n$-simplices of $D$.  A classical result asserts that this  chain
complex is acyclic, see for example \cite{We}, proof
of Theorem $8.3.8$. Therefore, we may equivalently consider the homology of
the normalised  simplicial chain complex $C^r(D)\doteq C(D)/C^d(D)$. This
chain complex is 
isomorphic  with the cellular chain complex of $|D|$, where $|D|$
is provided with its natural cell decomposition coming from the simplicial structure
of $D$. Recall that the geometric realisation $|D|$ of a
simplicial set  $D$ is a CW-complex with one $n$-cell for each
non-degenerate $n$-simplex of $D$. This can be used to prove  that the simplicial homology of a
simplicial set coincides with the cellular homology of its geometric
realisation.

In fact a stronger result holds. Namely, there exists an inclusion map
$i\colon C^r(D)\doteq C(D)/C^d(D) \to C(D)$ such that $p\circ i=\id$, with $i \circ
p$ being homotopic to the identity. Here $p\colon C(D) \to C^r(D)$ is the
projection map.  This is the Normalisation Theorem; see
\cite{ML}, VIII.6.

In this article, we will consider the unnormalised simplicial chain complex.

Note that if $M$ is an ordered  simplicial complex, then the simplicial chain complex of
$D_M$, the simplicial set made from $M$, coincides with  its usual
definition; see  for example \cite{Mau}, $4.3$. That book discusses both the
normalised and unnormalised cases.

\subsubsection{Homology of Crossed Modules}
Let $\G$ be a crossed module. The classifying space $B(\G)$  of $\G$ is the
geometric realisation of the simplicial set $N(\G)$, the nerve of $\G$. Therefore,  it is natural to consider its simplicial homology.

 The simplicial structure of $B(\G)$ was described in \ref{Refer2}. Let us
unpack the structure of the simplicial chain complex of $B(\G)$ for low
dimensions.  
Suppose that $\G=(G,E,\d,\t)$ is a reduced crossed module. 
For any $n\in \N$,  the simplicial chain group $C_n(\G)\doteq C_n(N(\G))$ is the free
abelian group on the set of all $\G$-colourings of $|\D(n)|$, with the obvious
boundary maps. In particular: 
\begin{enumerate}
\item  $C_0(\G)=\Z$,
\item  $C_1(\G)$ is the free $\Z$-module on the symbols $s(X),X \in G$,\item
  $C_2(\G)$  the  free $\Z$-module on the symbols $s(X,Y,e)$, where $X,Y \in G$ and $e  \in E$,
\item $C_3(\G)$ is the free $\Z$-module  on  the symbols $s(X,Y,Z,e,f,g,h)$, where $X,Y,Z \in G$ and $e,f,g,h
  \in E$ must verify $ef=(X \t g) h$.
\end{enumerate} 
Moreover we have:
\begin{align*}
 \d(s(X))&=0,\\
 \d(s(X,Y,e))&=s(X)+s(Y)-s\left (\d(e)^{-1}XY\right),
\end{align*}
where $X,Y \in G$ and $e \in E$; and also:
\begin{align*}
\d(s(X,Y,Z,e,f,g,h))&=s(Y,Z,g)-s(\d(e)^{-1}XY,Z,f)\\&\quad\quad\quad\quad\quad\quad\quad+s(X,\d(g)^{-1}YZ,h) -s(X,Y,e),
\end{align*}
where $X,Y,Z \in G$ and $e,f,g,h \in E$ are such that $ef=(X \t g) h$.

\begin{figure}
\centerline{\relabelbox 
\epsfysize 9cm
\epsfbox{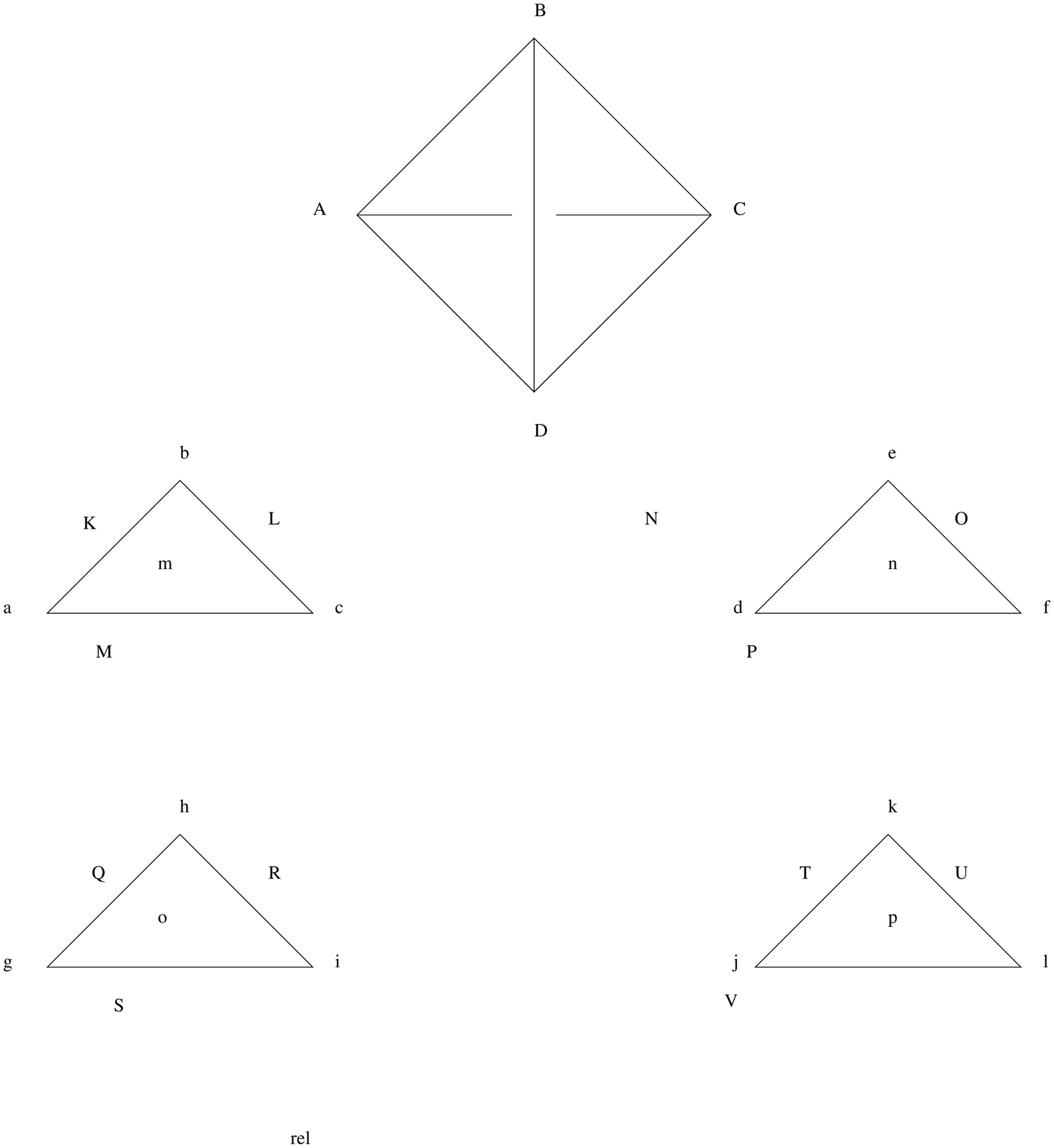}
\relabel {A}{$\s{0}$}
\relabel {B}{$\s{1}$}
\relabel {C}{$\s{2}$}
\relabel {D}{$\s{3}$}
\relabel {a}{$\s{0}$}
\relabel {b}{$\s{1}$}
\relabel {c}{$\s{2}$}
\relabel {d}{$\s{0}$}
\relabel {e}{$\s{2}$}
\relabel {f}{$\s{3}$}
\relabel {g}{$\s{1}$}
\relabel {h}{$\s{2}$}
\relabel {i}{$\s{3}$}
\relabel {j}{$\s{0}$}
\relabel {k}{$\s{1}$}
\relabel {l}{$\s{3}$}
\relabel {K}{$\s{X}$}
\relabel {L}{$\s{Y}$}
\relabel {M}{$\s{\d(e)^{-1}XY}$}
\relabel {N}{$\s{\d(e)^{-1}XY}$}
\relabel {O}{$\s{Z}$}
\relabel {P}{$\s{\d(f)^{-1}\d(e)^{-1}XYZ}$}
\relabel {Q}{$\s{Y}$}
\relabel {R}{$\s{Z}$}
\relabel {S}{$\s{\d(g)^{-1}YZ}$}
\relabel {T}{$\s{X}$}
\relabel {U}{$\s{\d(g)^{-1}YZ}$}
\relabel {V}{$\s{\d(h)^{-1}X\d(g)^{-1}YZ}$}
\relabel {rel}{$\s{ef=(X \t g)h}$}
\relabel {m}{$\s{e}$}
\relabel {n}{$\s{f}$}
\relabel {o}{$\s{g}$}
\relabel {p}{$\s{h}$}
\endrelabelbox }
\caption{\label{T1} The most general $\G$-colouring of $(01234)$: restriction
  to $(0123)$. }
\end{figure}

\begin{figure}
\centerline{\relabelbox 
\epsfysize 9cm
\epsfbox{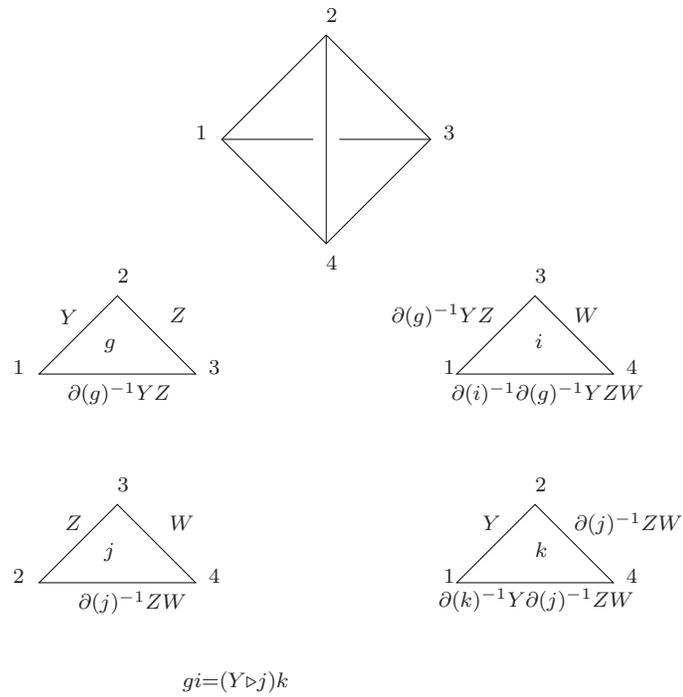}
\relabel {A}{$\s{1}$}
\relabel {B}{$\s{2}$}
\relabel {C}{$\s{3}$}
\relabel {D}{$\s{4}$}
\relabel {a}{$\s{1}$}
\relabel {b}{$\s{2}$}
\relabel {c}{$\s{3}$}
\relabel {d}{$\s{1}$}
\relabel {e}{$\s{3}$}
\relabel {f}{$\s{4}$}
\relabel {g}{$\s{2}$}
\relabel {h}{$\s{3}$}
\relabel {i}{$\s{4}$}
\relabel {j}{$\s{1}$}
\relabel {k}{$\s{2}$}
\relabel {l}{$\s{4}$}
\relabel {K}{$\s{Y}$}
\relabel {L}{$\s{Z}$}
\relabel {M}{$\s{\d(g)^{-1}YZ}$}
\relabel {N}{$\s{\d(g)^{-1}YZ}$}
\relabel {O}{$\s{W}$}
\relabel {P}{$\s{\d(i)^{-1}\d(g)^{-1}YZW}$}
\relabel {Q}{$\s{Z}$}
\relabel {R}{$\s{W}$}
\relabel {S}{$\s{\d(j)^{-1}ZW}$}
\relabel {T}{$\s{Y}$}
\relabel {U}{$\s{\d(j)^{-1}ZW}$}
\relabel {V}{$\s{\d(k)^{-1}Y\d(j)^{-1}ZW}$}
\relabel {rel}{$\s{gi=(Y \t j)k}$}
\relabel {m}{$\s{g}$}
\relabel {n}{$\s{i}$}
\relabel {o}{$\s{j}$}
\relabel {p}{$\s{k}$}
\endrelabelbox }
\caption{\label{T2} The most general $\G$-colouring of $(01234)$: restriction
  to $(1234)$.}
\end{figure}

\begin{figure}
\centerline{\relabelbox 
\epsfysize 9cm
\epsfbox{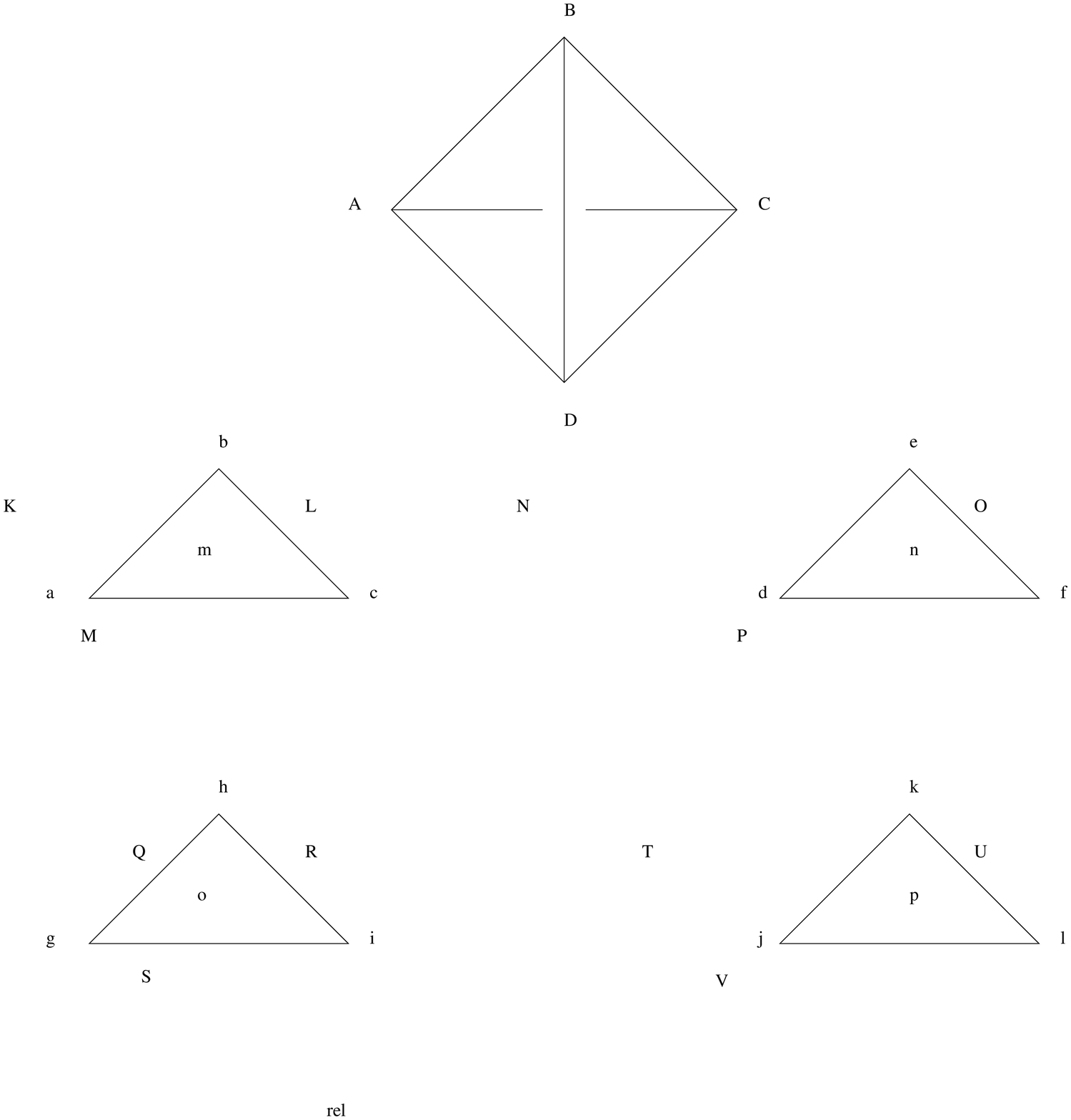}
\relabel {A}{$\s{0}$}
\relabel {B}{$\s{2}$}
\relabel {C}{$\s{3}$}
\relabel {D}{$\s{4}$}
\relabel {a}{$\s{0}$}
\relabel {b}{$\s{2}$}
\relabel {c}{$\s{3}$}
\relabel {d}{$\s{0}$}
\relabel {e}{$\s{3}$}
\relabel {f}{$\s{4}$}
\relabel {g}{$\s{2}$}
\relabel {h}{$\s{3}$}
\relabel {i}{$\s{4}$}
\relabel {j}{$\s{0}$}
\relabel {k}{$\s{2}$}
\relabel {l}{$\s{4}$}
\relabel {K}{$\s{\d(e)^{-1}XY}$}
\relabel {L}{$\s{Z}$}
\relabel {M}{$\s{\d(f)^{-1}\d(e)^{-1}XYZ}$}
\relabel {N}{$\s{\d(f)^{-1}\d(e)^{-1}XYZ}$}
\relabel {O}{$\s{W}$}
\relabel {P}{$\s{\d(m)^{-1}\d(f)^{-1}\d(e)^{-1}XYZW}$}
\relabel {Q}{$\s{Z}$}
\relabel {R}{$\s{W}$}
\relabel {S}{$\s{\d(j)^{-1}ZW}$}
\relabel {T}{$\s{\d(e)^{-1}XY}$}
\relabel {U}{$\s{\d(j)^{-1}ZW}$}
\relabel {V}{$\s{\d(n)^{-1}\d(e)^{-1}XY\d(j)^{-1}ZW}$}
\relabel {rel}{$\s{fm=e^{-1}(XY\t j) e n}$}
\relabel {m}{$\s{f}$}
\relabel {n}{$\s{m}$}
\relabel {o}{$\s{j}$}
\relabel {p}{$\s{n}$}
\endrelabelbox }
\caption{\label{T3} The most general $\G$-colouring of $(01234)$: restriction
  to $(0234)$. }
\end{figure}

\begin{figure}
\centerline{\relabelbox 
\epsfysize 9cm
\epsfbox{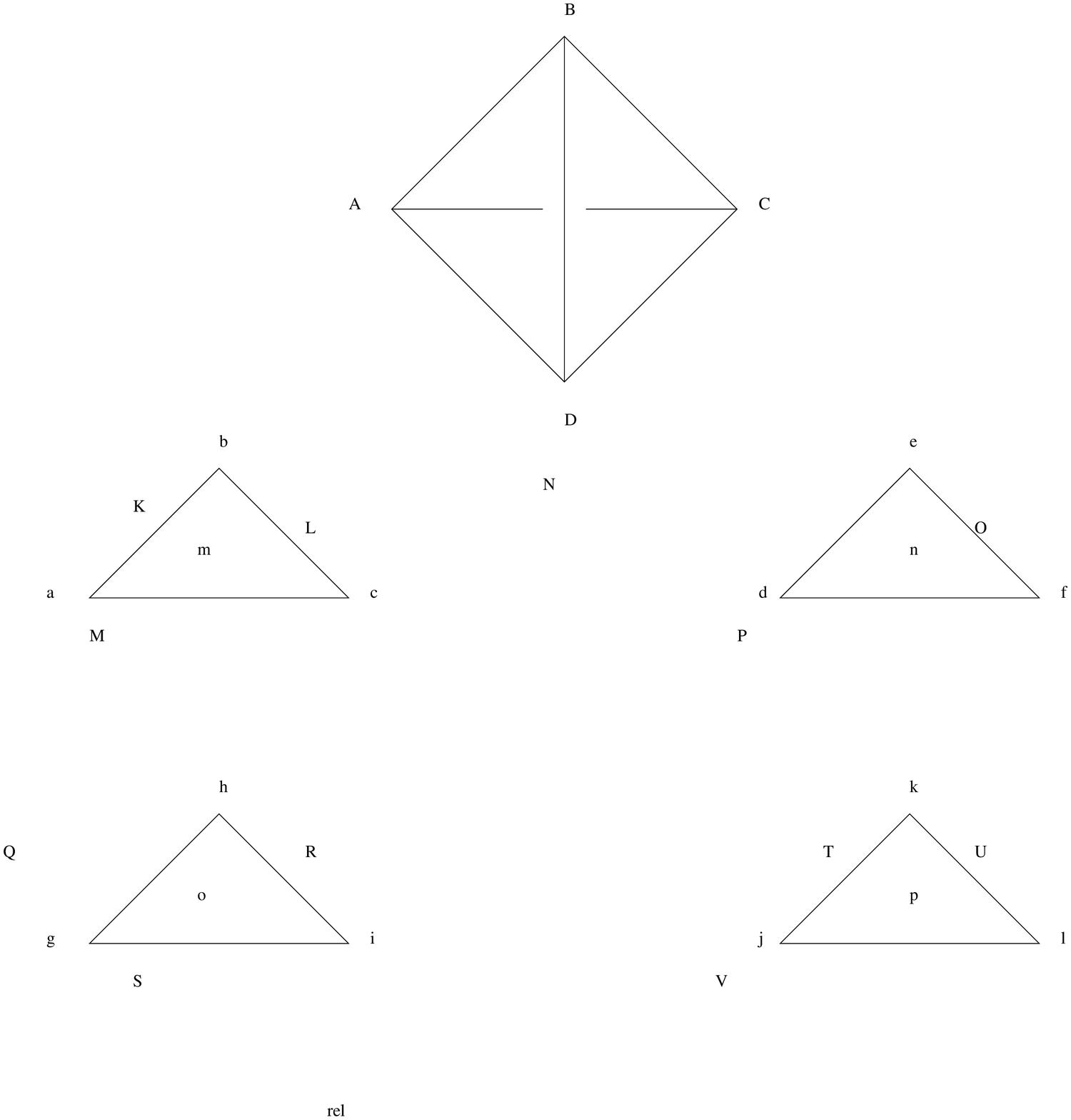}
\relabel {A}{$\s{0}$}
\relabel {B}{$\s{1}$}
\relabel {C}{$\s{3}$}
\relabel {D}{$\s{4}$}
\relabel {a}{$\s{0}$}
\relabel {b}{$\s{1}$}
\relabel {c}{$\s{3}$}
\relabel {d}{$\s{0}$}
\relabel {e}{$\s{3}$}
\relabel {f}{$\s{4}$}
\relabel {g}{$\s{1}$}
\relabel {h}{$\s{3}$}
\relabel {i}{$\s{4}$}
\relabel {j}{$\s{0}$}
\relabel {k}{$\s{1}$}
\relabel {l}{$\s{4}$}
\relabel {K}{$\s{X}$}
\relabel {L}{$\s{\d(g)^{-1}YZ}$}
\relabel {M}{$\s{\d(h)^{-1}X \d(g)^{-1}YZ}$}
\relabel {N}{$\s{\d(f)^{-1}\d(e)^{-1}XYZ}$}
\relabel {O}{$\s{W}$}
\relabel {P}{$\s{\d(m)^{-1}\d(f)^{-1}\d(e)^{-1}XYZW}$}
\relabel {Q}{$\s{\d(g)^{-1}YZ}$}
\relabel {R}{$\s{W}$}
\relabel {S}{$\s{\d(i)^{-1}\d(g)^{-1}YZW}$}
\relabel {T}{$\s{X}$}
\relabel {U}{$\s{\d(i)^{-1}\d(g)^{-1}YZW}$}
\relabel {V}{$\s{\d(p)^{-1}X\d(i)^{-1}\d(g)^{-1}YZW}$}
\relabel {rel}{$\s{hm=(X \t i)p}$}
\relabel {m}{$\s{h}$}
\relabel {n}{$\s{m}$}
\relabel {o}{$\s{i}$}
\relabel {p}{$\s{p}$}
\endrelabelbox }
\caption{\label{T4} The most general $\G$-colouring of $(01234)$: restriction
  to $(0134)$. }
\end{figure}

\begin{figure}
\centerline{\relabelbox 
\epsfysize 9cm
\epsfbox{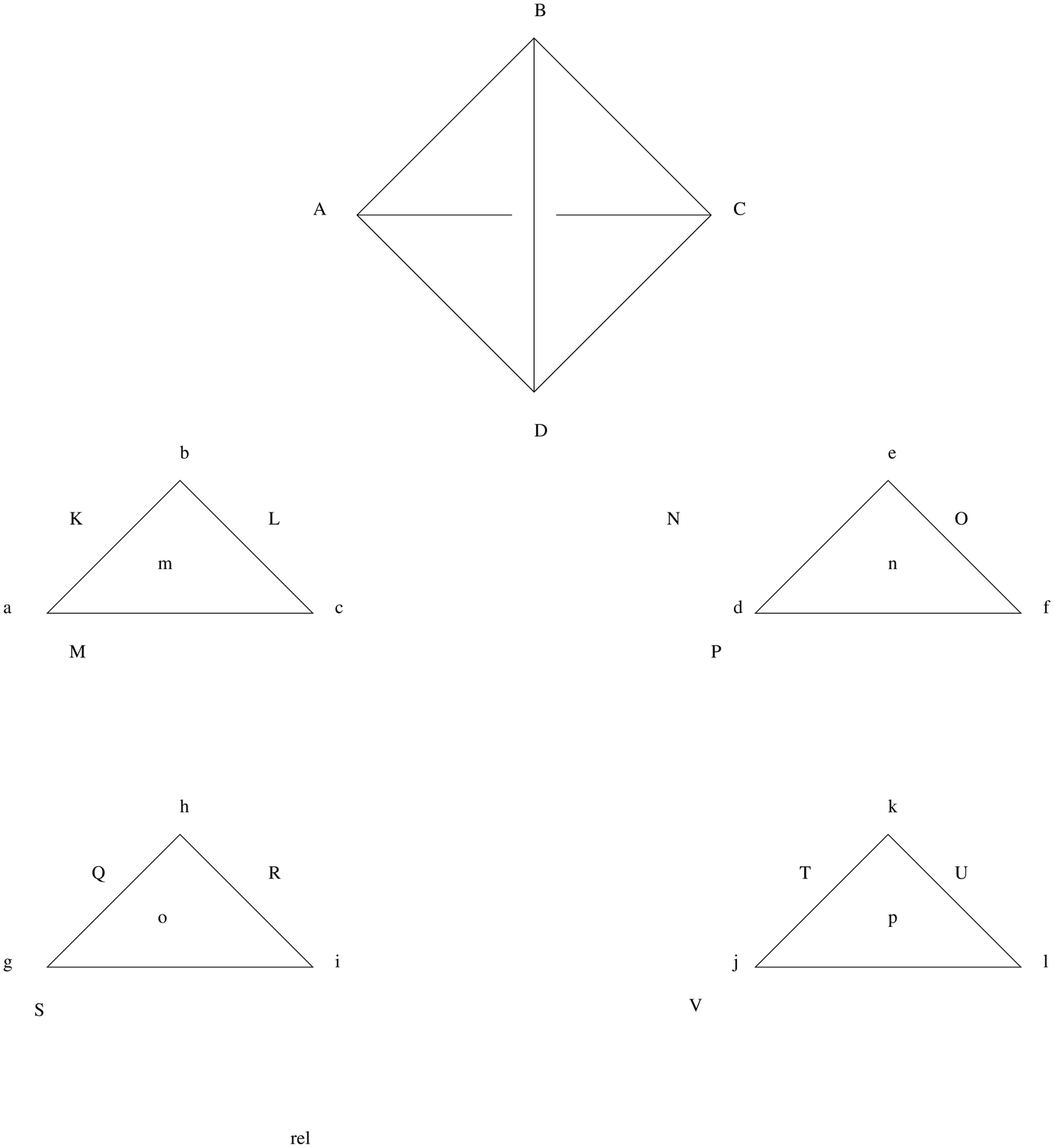}
\relabel {A}{$\s{0}$}
\relabel {B}{$\s{1}$}
\relabel {C}{$\s{2}$}
\relabel {D}{$\s{4}$}
\relabel {a}{$\s{0}$}
\relabel {b}{$\s{1}$}
\relabel {c}{$\s{2}$}
\relabel {d}{$\s{0}$}
\relabel {e}{$\s{2}$}
\relabel {f}{$\s{4}$}
\relabel {g}{$\s{1}$}
\relabel {h}{$\s{2}$}
\relabel {i}{$\s{4}$}
\relabel {j}{$\s{0}$}
\relabel {k}{$\s{1}$}
\relabel {l}{$\s{4}$}
\relabel {K}{$\s{X}$}
\relabel {L}{$\s{Y}$}
\relabel {M}{$\s{\d(e)^{-1}XY}$}
\relabel {N}{$\s{\d(e)^{-1}XY}$}
\relabel {O}{$\s{\d(j)^{-1}ZW}$}
\relabel {P}{$\s{\d(n)^{-1}\d(e)^{-1}XY\d(j)^{-1}ZW}$}
\relabel {Q}{$\s{Y}$}
\relabel {R}{$\s{\d(j)^{-1}ZW}$}
\relabel {S}{$\s{\d(k)^{-1}Y\d(j)^{-1}ZW}$}
\relabel {T}{$\s{X}$}
\relabel {U}{$\s{\d(i)^{-1}\d(g)^{-1}YZW}$}
\relabel {V}{$\s{\d(p)^{-1}X\d(i)^{-1}\d(g)^{-1}YZW}$}
\relabel {rel}{$\s{en=(X \t k)p}$}
\relabel {m}{$\s{e}$}
\relabel {n}{$\s{n}$}
\relabel {o}{$\s{k}$}
\relabel {p}{$\s{p}$}
\endrelabelbox }
\caption{\label{T5} The most general $\G$-colouring of $(01234)$: restriction
  to $(0124)$. }
\end{figure}

The determination of $C_4(\G)$ is a bit more complicated. In figures  \ref{T1}
to \ref{T5}, we display the most general $\G$-colouring of the $4$-simplex
$(01234)$. It depends on the variables $X,Y,Z,W \in G$ and
$e,f,g,h,i,j,k,m,n,p \in E$, which must satisfy the conditions  shown in figures  \ref{T1}
to \ref{T5}. Namely: $
 ef=(X\t g)h,$ $       gi=(Y \t j) k,$ $  fm=e^{-1}(XY \t j) e n,$ $  hm=(X \t i)p $ and $en=(X \t k)p$. Note that the last relation follows from all the others. From these relations,  it is
not difficult  to conclude that the colourings of
all $1$-simplices of figures  \ref{T1}
to \ref{T5} are coherent.
 The associated simplicial 4-chains   are denoted by
$s(X,Y,Z,W,e,f,g,h,i,j,k,m,n,p)$. The determination of  the boundary map
$\d\colon   C_4(\G) \to C_3(\G)$  is an easy task.

\subsubsection{Explicit Description of the Low Dimensional Coboundary Maps}
Let $\G=(G,E,\d,\t)$ be a reduced crossed module. We consider the
$U(1)$-cohomology of $\G$.  This is a very particular case of the construction
in \cite{Pa}. There was considered the general case of cohomology with coefficients in any $\pi_1(B(\G))$-module. From the discussion above, the group $C^3(\G)$ of $3$-cochains of $\G$ is given by all maps
$\w\colon    G^3 \times E^4 \to U(1)$ which verify:
\[\w(X,Y,Z,e,f,g,h) \neq 1 \implies ef=(X \t g) h.\]
Analogously, the group $C^4(\G)$ of $4$-cochains is given by all maps 
$G^4 \times E^{10} \to U(1)$ which verify:
\begin{equation*}
\w(X,Y,Z,W,e,f,g,h,i,j,k,m,n,p)\neq 1 \implies 
\left\{ \begin{CD} ef=(X\t g)h, \\
                   gi=(Y \t j) k,\\
                   fm=e^{-1}(XY \t j) e n, \\
                   hm=(X \t i)p,\\
                   en=(X \t k)p.
         \end{CD} \right.
\end{equation*}
Recall  that the last relation is a consequence of the others. The group
$C^2(\G)$ of $2$-cochains is simply given by all maps $\w\colon   G^2 \times E \to
U(1)$.   The following results follow trivially.

\begin{Proposition}\label{oneC}
Let $\w \in C^1(\G)$ be a 1-cochain. Then:
\[d{\w}(X,Y,e)=\w(X)+\w(Y)-\w(\d(e)^{-1} XY), \]
for all $ X,Y \in G$ and all $ e \in E.$ 
\end{Proposition} 

\begin{Proposition}\label{twoC}
Let $\w\in C^2(\G)$ be a 2-cochain. Then
\begin{multline*}
d\w(X,Y,Z,e,f,g,h)=\\
   \quad \quad \w(Y,Z,g)-\w(\d(e)^{-1}XY,Z,f)+\w(X,\d(g)^{-1}YZ,h)-\w(X,Y,e),
\end{multline*}
whenever $ef=(X\t g)h$ and equals  $1\in \mathbb{C}$ otherwise. 
\end{Proposition}
\begin{Proposition}\label{threeC}
Let $\w\in C^3(\G,A)$ be a 3-cochain. We have:
\begin{multline*}
d\w(X,Y,Z,W,e,f,g,h,i,j,k,m,n,p)=\\
\w(X,Y,Z,e,f,g,h)+\w(Y,Z,W,g,i,j,k)
+\w(X,\d(g)^{-1}YZ,W,h,m,i,p)\\-\w(\d(e)^{-1}XY,Z,W,f,m,j,n)
-\w(X,Y,\d(j)^ {-1}ZW,e,n,k,p),
\end{multline*}
whenever
$
 ef=(X\t g)h,$ $       gi=(Y \t j) k,$ $  fm=e^{-1}(XY \t j) e n,$ $  hm=(X \t i)p $ and $en=(X \t k)p,$
and equals $1\in \mathbb{C}$ otherwise.
\end{Proposition}

\subsection{A Homotopy Invariant of 3-Manifolds}
We will restrict our discussion to the $3$-dimensional case. However, it is clear that the results that we obtain will still
hold, with the obvious adaptations,  for any dimension  $n\in \N$, and can be
extended to handle crossed complexes {in the same way}. The $n$-dimensional analogues of propositions \ref{oneC} to \ref{threeC}, as well as their extension to crossed complexes, important for calculational purposes, require, however, more laborious calculations.

 Let  $M$ be a 3-dimensional oriented triangulated closed piecewise linear manifold.  The orientation class
 $o_M\in  H_3(M)$ of $M$ chosen can be specified  by an assignment of a
 total order  to each non-degenerate tetrahedron of $M$. These total orders are  defined up to even
 permutations. They define an orientation on each tetrahedron of $M$.  Therefore it is  required that 
 if two non-degenerate
 tetrahedra share a non-degenerate face then the orientations induced on their common face
 should be opposite.

  As usual, we suppose that we are provided with a total order on the set of all  vertices of
 $M$. Consequently, each non-degenerate $3$-simplex $K$ of $M$ can be uniquely represented as $K=(abcd)$ where
 $a<b<c<d$. If $K$ is a non-degenerate tetrahedron, we say that $r(K)$ is $-1$ or $1$ according to whether the total
 order induced on the vertices of $K$ differs from the one determined by the
 orientation of $M$ by an odd or an even permutation. {In other words, $r(K)$ is $1$ or $-1$ depending on whether the  orientation on $K$ induced by the total order on the set of vertices of $M$ coincides or not with the orientation of $K$ determined by the orientation of $M$.}

The orientation class $o_M$ of $M$, living in the (normalised or unnormalised)
simplicial  homology group
$H_3(M)$ of $M$ is therefore: 
\[o_M=\sum_{\textrm{3-simplices }(abcd)}r(abcd) {(abcd)},\]
where the sum is extended to the non-degenerate 3-simplices, only.

Let $\G=(G,E,\d,\t)$ be a finite reduced crossed module. Choose a 3-dimensional cocycle $\w$ representing some cohomology class in $ H^3(\G)$. For any $\G$-colouring
$\bc$ of $M$ define the $U(1)$-valued ``action'':
\begin{multline}
\S(\bc,\w)\\=\prod_{\textrm{3-simplices }(abcd)}
\w\big(\bc(ab),\bc(bc),\bc(cd),\bc(abc),
\bc(acd),\bc(bcd),\bc(abd)\big)^{r(abcd)}, \end{multline}
where the product is extended to the non-degenerate 3-simplices of $M$, only.

Recall that the set of $3$-simplices of the classifying space $B(\G)$ of $\G$ is in one-to-one correspondence with the set of $\G$-colourings of the standard geometric 3-simplex $|\D(3)|$.  The group of 3-dimensional simplicial cochains of $B(\G)$ is given by all assignments of an element of $U(1)$ to each $\G$-colouring $\bc$ of $|\D(3)|$. Given a non-degenerate 3-simplex $(abcd)$ of $M$, the quantity  $\w\big(\bc(ab),\bc(bc),\bc(cd),\bc(abc), \bc(acd),\bc(bcd),\bc(abd)\big)$ is by definition, and under the identification above, exactly $\w(\bc_{|(abcd)})$, where $\bc_{|(abcd)}$ is the restriction of the $\G$-colouring $\bc$ of $M$ to the tetrahedron $(abcd)$.

\begin{Theorem}\label{main2}
Let $M$ be a 3-dimensional closed oriented triangulated piecewise linear  manifold, with a total order
on   its set of  vertices. Let $n_0$ and $n_1$ be, respectively, the number
of vertices and edges of $M$. Let also $\G=(G,E,\d,\t)$ be a finite reduced crossed
module, and let ${\w} \in H^3(\G)$ be a 3-dimensional cohomology class of $\G$. 
The quantity:
\[I_\G(M,\w)=\frac{\#E^{n_0}}{\#G^{n_0}\#E^{n_1}}\sum_{\G\textrm{-colourings } \bc } \S(\bc,\w)\]
is a homotopy invariant of $M$, and therefore, in particular, it is independent of 
the ordered triangulation of $M$ chosen. In fact, let $o_M \in H_3(M) $ be the orientation
class of $M$. We have:
\[
I_\G(M,\w)=\sum_{g \in [M,B(\G)]}\frac{ \#\pi_2(\TOP(M,B(\G)),g)}{
  \#\pi_1(\TOP(M,B(\G)),g)}  \left <o_M,g^*(\w)        \right>.
\]
Here $[M,B(\G)]=\pi_0(\TOP(M,B(\G)))$ denotes the set of homotopy classes of
maps $M \to B(\G)$.
\end{Theorem}
Note that since $N(\G)$ is Kan,  the Simplicial Approximation Theorem
guarantees that  any map $f\colon    M \to B(\G)$ is homotopic to the geometric
realisation of  a simplicial map $T_M \to N(\G)$,
defined up to simplicial homotopy. Here $T_M$ is the simplicial
set defined from the triangulation of $M$. In particular $f^*(\w)$ is well defined in the simplicial category for any continuous map $f\colon    M \to B(\G)$.
The Simplicial Approximation Theorem (for simplicial sets) is proved for example in \cite{S}.  Note also the Normalisation Theorem stated in \ref{homKan}.

The proof of Theorem \ref{main2} is analogous to the proof of Theorem
\ref{main1}. The main lemma which we will use for its proof is the following.

\begin{Lemma}\label{C}
Let $f \in \CRS_0(\Pi(M),\G)$ be a morphism $\Pi(M) \to \G$. Therefore, by Proposition \ref{Refer1} we can
associate a $\G$-colouring $\bc^f$ of $M$ to it.  We have:
\[\left <o_M,\n(f)^*({\w})\right>=\S\left (\bc^f,\w\right).\]
Note that $\n(f)$ is the realisation of a simplicial map. In fact it is the geometric
realisation of $F(f)$; see theorems \ref{Refer4} and \ref{weak}.
\end{Lemma}
In particular, from theorems \ref{weak} and \ref{fundamental} and  subsequent comments, it follows that
$S\left (\bc^f,\w\right)$  depends only on the homotopy class of $f\colon   
\Pi(M) \to \G$.  This also proves that the action $S({\bf c}^f,\w)$
does not depend on the  cocycle representing the cohomology class $\w\in H^3(\G)$.

\begin{Proof}  {\bf (Lemma \ref{C})}
Recall the notation introduced in \ref{Refer2}.  The $\G$-colouring $\bc^f
$ of
$M$ restricts to a $\G$-colouring $\bc^f_{|K}$ of $K$, for each non-degenerate simplex $K$  of $M$.
We have: 
\begin{align*}
&\left <o_M,\n(f)^*\w\right>\\
&=\left <o_M,F(f)^*\w\right>, \textrm{by definition (cf. Theorem \ref{weak})}\\
&=\left <F(f)_*o_M,\w\right)\\
&= \left < \sum_{\textrm{3-simplices }(abcd)} r(abcd) F(f)(abcd),{\w}\right>\\
&= \prod_{\textrm{3-simplices }(abcd)} \left
  <F(f)(abcd),{\w}\right>^{r(abcd)}\\
&= \prod_{\textrm{3-simplices }(abcd)} \left  <f_{|\Pi(abcd)},{\w}\right>^{r(abcd)},\textrm{ by Theorem \ref{Refer4} and Remark \ref{REFER}}\\
&= \prod_{\textrm{3-simplices }(abcd)}\left <{\bc^f}_{|(abcd)},{\w}\right>^{r(abcd)},
\textrm { by Proposition \ref{Refer1}}\\
&=\prod_{\textrm{3-simplices }(abcd)}\hskip -6mm\w\big(\bc^f(ab),\bc^f(bc),\bc^f(cd),\bc^f(abc), \bc^f(acd),\bc^f(bcd),\bc^f(abd)\big)^{r(abcd)}\\
&=\S\left (\bc^f,\w \right ).
\end{align*}
The second to last step follows by definition. Note that the sum and the
products are to be extended to the non-degenerate 3-simplices of $M$, only.
\end{Proof}

We now prove Theorem \ref{main2}. 

\begin{Proof} {\bf (Theorem \ref{main2}) }
We maintain the notation we used in  the proof of Theorem \ref{main1}. We have:
\begin{align*}
\sum_{g \in [M,B(\G)]}
&\frac{\#\pi_2(\TOP(M,B(\G)),g)}{\#\pi_1(\TOP(M,B(\G)),g)}\left <o_M,g^*\w\right>
\\&=\sum_{f \in \CRS_0(\Pi(M),\G)}
\frac{\#\pi_2(\TOP(M,B(\G)),\n(f))}{\#\pi_1(\TOP(M,B(\G)),\n(f)) \#[f]}
\left <o_M,\n(f)^*\w \right > 
\\&=\sum_{f \in \CRS_0(\Pi(M),\G)}
\frac{\#\pi_2(\TOP(M,B(\G)),\n(f))}{\#\pi_1(\TOP(M,B(\G)),\n(f)) \#[f]} \S\left (\bc^f,\w\right).
\end{align*}
The same calculation as in the proof of theorem \ref{main1} finishes the proof. Note that we are implicitly using the fact that if $f$ and $f'$ belong to the same connected component in the groupoid $\CRS_1(\Pi(M),\G)$ then it follows that 
$\n(f)$ is homotopic to $\n(f')$ and $\S\left (c^f,\w\right )=\S\left (c^{f'},\w\right)$.
\end{Proof}

As we referred to before, this theorem can be  extended in the obvious way to closed $n$-manifolds, with $n$ arbitrary, and
cohomology classes of crossed complexes. Compare with Theorem \ref{Main3}. It would be interesting to relate our construction with M. Mackaay's  work appearing in \cite{Mk2}. Conjecturally, this last should be related  to the 4-manifold invariant obtained from $4$-dimensional cohomology classes of crossed complexes of length $3$. Finding the precise link  forces the determination of all the relations verified by  4-dimensional crossed complex cocycles, which itself requires elaborate calculations. We will consider these issues in a subsequent publication.

\section*{Acknowledgements}
The authors would like to thank  Ronnie Brown, Gustavo Granja, Roger Picken, Marco Mackaay  and Jim Stasheff  for useful comments.

\end{document}